\documentclass[12pt,twoside,english,british]{article}
\usepackage[T1]{fontenc}
\usepackage[utf8]{inputenc}
\usepackage[a4paper]{geometry}
\usepackage{color}
\usepackage{amsmath}
\usepackage{amsthm}
\usepackage{amssymb}
\usepackage{setspace}
\usepackage{esint}

\makeatletter
\numberwithin{equation}{section}
\numberwithin{figure}{section}
\theoremstyle{plain}
\newtheorem{thm}{\protect\theoremname}[section]
\theoremstyle{definition}
\newtheorem{defn}[thm]{\protect\definitionname}
\theoremstyle{plain}
\newtheorem{lem}[thm]{\protect\lemmaname}
\theoremstyle{plain}
\newtheorem{prop}[thm]{\protect\propositionname}
\newcommand{\lyxaddress}[1]{
	\par {\raggedright #1
	\vspace{1.4em}
	\noindent\par}
}

\usepackage[T1]{fontenc}
\usepackage[utf8]{inputenc}
\usepackage[a4paper]{geometry}
\usepackage{color}
\usepackage{dsfont}
\usepackage{amsmath}
\usepackage{amsthm}
\usepackage{amssymb}
\usepackage{setspace}
\usepackage{esint}

\makeatletter
\numberwithin{equation}{section}
\numberwithin{figure}{section}
\theoremstyle{plain}

\@ifundefined{date}{}{\date{}}
\usepackage{babel}
\usepackage{latexsym}
\usepackage{amsthm}
\usepackage{esint}
\usepackage{eucal}
\usepackage{epstopdf}
\usepackage{graphicx}
\usepackage{bigints}
\theoremstyle{plain}

\newtheoremstyle{boldremark}
    {\dimexpr\topsep/2\relax} 
    {\dimexpr\topsep/2\relax} 
    {}          
    {}          
    {\bfseries} 
    {.}         
    {.5em}      
    {}          

\theoremstyle{boldremark}
\newtheorem{brem} [thm] {Remark} 

\usepackage{fancyhdr}
\usepackage{blindtext}
\usepackage{titlesec}

\titleformat{name=\chapter}[display]
{\normalfont\Huge\itshape}
{\titlerule[1pt]\vspace{-33pt}\filleft%
  \parbox[t]{6em}{%
    \raggedleft%
    \rule{\linewidth}{0.5ex}\newline%
    \chaptertitlename\ \thechapter%
  }%
}
{4pc}
{\normalfont\upshape\bfseries\Huge}
\allowdisplaybreaks

\makeatother

\usepackage{babel}
\addto\captionsbritish{\renewcommand{\definitionname}{Definition}}
\addto\captionsbritish{\renewcommand{\lemmaname}{Lemma}}
\addto\captionsbritish{\renewcommand{\theoremname}{Theorem}}
\addto\captionsenglish{\renewcommand{\definitionname}{Definition}}
\addto\captionsenglish{\renewcommand{\lemmaname}{Lemma}}
\addto\captionsenglish{\renewcommand{\theoremname}{Theorem}}
\providecommand{\definitionname}{Definition}
\providecommand{\lemmaname}{Lemma}
\providecommand{\theoremname}{Theorem}

\usepackage[colorlinks,pdfpagelabels,pdfstartview = FitH,bookmarksopen
= true,bookmarksnumbered = true,linkcolor = red,plainpages =
false,hypertexnames = false,citecolor = blue,pagebackref=false]{hyperref}

\makeatother

\usepackage{babel}
\addto\captionsbritish{\renewcommand{\definitionname}{Definition}}
\addto\captionsbritish{\renewcommand{\lemmaname}{Lemma}}
\addto\captionsbritish{\renewcommand{\propositionname}{Proposition}}
\addto\captionsbritish{\renewcommand{\theoremname}{Theorem}}
\addto\captionsenglish{\renewcommand{\definitionname}{Definition}}
\addto\captionsenglish{\renewcommand{\lemmaname}{Lemma}}
\addto\captionsenglish{\renewcommand{\propositionname}{Proposition}}
\addto\captionsenglish{\renewcommand{\theoremname}{Theorem}}
\providecommand{\definitionname}{Definition}
\providecommand{\lemmaname}{Lemma}
\providecommand{\propositionname}{Proposition}
\providecommand{\theoremname}{Theorem}

\begin{document}
\title{\textbf{Sharp Sobolev regularity for widely\\degenerate parabolic
equations}}
\author{Pasquale Ambrosio \thanks{Dipartimento di Matematica, Università di Bologna, Piazza di Porta
S. Donato 5, 40126 Bologna, Italy.\textit{ E-mail address}: pasquale.ambrosio@unibo.it}}
\date{November 22, 2024}

\maketitle

\begin{abstract}
\noindent We consider local weak solutions to the widely degenerate
parabolic PDE 
\[
\partial_{t}u-\mathrm{div}\left((\vert Du\vert-\lambda)_{+}^{p-1}\frac{Du}{\vert Du\vert}\right)=f\,\,\,\,\,\,\,\,\mathrm{in}\,\,\,\Omega_{T}=\Omega\times(0,T),
\]
where $p\geq2$, $\Omega$ is a bounded domain in $\mathbb{R}^{n}$
for $n\geq2$, $\lambda$ is a non-negative constant and $\left(\,\cdot\,\right)_{+}$
stands for the positive part. Assuming that the datum $f$ belongs
to a suitable Lebesgue-Besov parabolic space when $p>2$ and that
$f\in L_{loc}^{2}(\Omega_{T})$ if $p=2$, we prove the Sobolev spatial
regularity of a\textit{ novel} nonlinear function of the spatial gradient
of the weak solutions. This result, in turn, implies the existence
of the weak time derivative for the solutions of the evolutionary
$p$-Poisson equation. The main novelty here is that $f$ only has
a Besov or Lebesgue spatial regularity, unlike the previous work \cite{AmPa},
where $f$ was assumed to possess a Sobolev spatial regularity of
integer order. We emphasize that the results obtained here can be
considered, on the one hand, as the parabolic analog of some elliptic
results established in \cite{AmGrPa}, and on the other hand as the
extension to a strongly degenerate setting of some known results for
less degenerate parabolic equations.\vspace{0.2cm}
 
\end{abstract}


\section{Introduction }

\noindent $\hspace*{1em}$In this paper, we are interested in local
weak solutions $u:\Omega_{T}\rightarrow\mathbb{R}$ to the evolutionary
equation 
\begin{equation}
\partial_{t}u-\mathrm{div}\left((\vert Du\vert-\lambda)_{+}^{p-1}\frac{Du}{\vert Du\vert}\right)=f\,\,\,\,\,\,\,\,\mathrm{in}\,\,\,\Omega_{T}=\Omega\times(0,T),\label{eq:1}
\end{equation}
where $\Omega$ is a bounded domain in $\mathbb{R}^{n}$ ($n\geq2$),
$T>0$, $\lambda\geq0$ and $\left(\,\cdot\,\right)_{+}$ stands for
the positive part. Our analysis is confined to the super-quadratic
case, i.e. in this work we consider $p\geq2$.\\
$\hspace*{1em}$A motivation for studying PDEs of the type (\ref{eq:1})
can be found in gas filtration problems taking into account the initial
pressure gradient (see \cite{AmCuDe}, \cite{AmPa}, \cite{BDGP0}
and the references therein). 

\noindent $\hspace*{1em}$The main feature of the above equation is
that it exhibits a strong degeneracy, coming from the fact that its
modulus of ellipticity vanishes in the region $\{\vert Du\vert\leq\lambda\}$,
and hence its principal part behaves like a $p$-Laplace operator
only for large values of $\vert Du\vert$.\\
$\hspace*{1em}$As our main result, here we establish the Sobolev
spatial regularity of a nonlinear function of the spatial gradient
$Du$ of the weak solutions to equation (\ref{eq:1}), by assuming
that the datum $f$ belongs to a suitable Lebesgue-Besov parabolic
space when $p>2$ (see Theorem \ref{thm:theo1} below) and that $f\in L_{loc}^{2}(\Omega_{T})$
if $p=2$ (see Theorem \ref{thm:nuovo}). These results, in turn,
imply the Sobolev time regularity of the weak solutions to the evolutionary
$p$-Poisson equation, under the same hypotheses on the function $f$
(cf. Theorem \ref{thm:theo2}, where we only address the case $p>2$,
since the Sobolev time regularity is well known for the heat equation
with source term in $L_{loc}^{2}(\Omega_{T})$). These issues have
been widely investigated, as one can see, for example, in \cite{AmPa,Duzaar,GPdN,GiPaSc,Sche}.
In particular, establishing the Sobolev regularity of the solutions
with respect to time, once the higher differentiability in space has
been obtained, is a quite standard fact: see, for instance, \cite{Lind1}-\cite{Lind3}.\\
$\hspace*{1em}$Before specifying in detail the assumption on $f$
in the case $p>2$, we wish to discuss some results already available
in the literature. A common aspect of nonlinear parabolic problems
with growth rate $p\geq2$ is that the higher differentiability is
proven for a nonlinear function of the gradient that takes into account
the growth of the principal part of the equation. Indeed, already
for the parabolic $p$-Laplace equation (which is obtained from (\ref{eq:1})
by setting $\lambda=f=0$), the higher differentiability is established
for the function 
\[
V(Du):=\vert Du\vert^{\frac{p-2}{2}}Du\,,
\]
as one can see, for example, in \cite{Duzaar}. In case of widely
degenerate problems, this phenomenon persists and, both for elliptic
and parabolic equations, higher differentiability results hold true
for the function $H_{p/2}(Du):=(\vert Du\vert-\lambda)_{+}^{p/2}\frac{Du}{\left|Du\right|}$
(see \cite{Amb1,Amb2,AmPa,Br}). However, this function does not provide
any information about the second-order regularity of the solutions
in the set where the equation becomes degenerate. Actually, since
every $\lambda$-Lipschitz function is a solution of the homogeneous
elliptic equation 
\[
\mathrm{div}\left((\vert Du\vert-\lambda)_{+}^{p-1}\frac{Du}{\vert Du\vert}\right)=0\,,
\]
no more than Lipschitz regularity can be expected for the solutions:
in this regard, see \cite{BrLIP} for the elliptic counterpart of
equation (\ref{eq:1}) and \cite{AmBa} for a very recent generalization
to strongly degenerate parabolic systems.

\noindent $\hspace*{1em}$In addition, it is well known that, for
both the stationary and the evolutionary $p$-Poisson equation, a
Sobolev (spatial) regularity is required for the datum $f$ in order
to get the Sobolev (spatial) regularity of the function $V(Du)$,
which coincides with $H_{p/2}(Du)$ for $\lambda=0$ (see, for example,
\cite{BraSan} for the stationary $p$-Poisson equation and \cite{AmPa}
for the evolutionary one). In particular, in the elliptic setting,
the optimal assumption on the datum $f$ to obtain the local $W^{1,2}$-regularity
of $V(Du)$ has been determined in \cite[Theorem 1.1]{BraSan} as
a fractional Sobolev regularity suitably related to the growth exponent
$p>2$. Recently, this result has somewhat been improved and extended
to the elliptic version of (\ref{eq:1}) in \cite{AmGrPa}: there,
for $p>2$, $\lambda>0$ and $\alpha\geq\frac{p+1}{2(p-1)}$, the
authors consider the function $\mathcal{V}_{\alpha,\lambda}:\mathbb{R}^{n}\rightarrow\mathbb{R}^{n}$
defined by 
\begin{equation}
\mathcal{V}_{\alpha,\lambda}(\xi)\equiv\mathcal{V}_{\lambda}(\xi):=\begin{cases}
\begin{array}{cc}
\mathcal{G}_{\alpha,\lambda}((\vert\xi\vert-\lambda)_{+})\,\frac{\xi}{\left|\xi\right|} & \,\,\,\mathrm{if}\,\,\,\vert\xi\vert>\lambda,\\
0 & \,\,\,\mathrm{if}\,\,\,\vert\xi\vert\leq\lambda,
\end{array}\end{cases}\label{eq:Vfun}
\end{equation}
with
\begin{equation}
\mathcal{G}_{\alpha,\lambda}(t):=\int_{0}^{t}\frac{\omega^{\frac{p-1+2\alpha}{2}}}{(\omega+\lambda)^{\frac{1+2\alpha}{2}}}\,d\omega\,\,\,\,\,\,\,\,\,\,\mathrm{for}\,\,t\geq0,\label{eq:Gfun}
\end{equation}
and establish that $\mathcal{V}_{\alpha,\lambda}(Du)\in W_{loc}^{1,2}(\Omega,\mathbb{R}^{n})$
by assuming 
\[
f\in B_{p',1,loc}^{\frac{p-2}{p}}(\Omega)\,\,\,\,\,\,\,\,\,\,\mathrm{with}\,\,\,p':=\frac{p}{p-1}\,,
\]
which is weaker than the assumption adopted in \cite{BraSan}, due
to Lemma \ref{lem:emb} below.\\
$\hspace*{1em}$Here, our main results follow in the footsteps of
the ones mentioned above, with the aim of extending the elliptic results
contained in \cite{AmGrPa} to the parabolic setting. The first theorem
of this work reads as follows. We refer to Sections \ref{sec:prelim}
and \ref{sec:Functional spaces} for further notations and definitions. 
\begin{thm}
\noindent \label{thm:theo1} Let $n\geq2$, $p>2$, $\lambda\geq0$
and $f\in L_{loc}^{p'}\left(0,T;B_{p',1,loc}^{\frac{p-2}{p}}(\Omega)\right)$.
Moreover, let 
\begin{equation}
\alpha=\begin{cases}
\begin{array}{cc}
0 & \textit{if}\,\,\,\lambda=0,\\
\textit{any\,\,value\,\,in}\,\,\left[\frac{p+1}{2(p-1)},\infty\right) & \textit{if}\,\,\,\lambda>0,
\end{array}\end{cases}\label{eq:alpha}
\end{equation}
and assume that 
\[
u\in C^{0}\left(0,T;L^{2}(\Omega)\right)\cap L_{loc}^{p}\left(0,T;W_{loc}^{1,p}(\Omega)\right)
\]
is a local weak solution of equation $\mathrm{(\ref{eq:1})}$. Then
\[
\mathcal{V}_{\alpha,\lambda}(Du)\,\in\,L_{loc}^{2}\left(0,T;W_{loc}^{1,2}(\Omega,\mathbb{R}^{n})\right).
\]
Furthermore, for any parabolic cylinder $Q_{r}(z_{0})\subset Q_{\rho}(z_{0})\subset Q_{R}(z_{0})\Subset\Omega_{T}$
we have\begin{align}\label{eq:apriori1}
\int_{Q_{r/2}(z_{0})}\left|D_{x}\mathcal{V}_{\alpha,\lambda}(Du)\right|^{2}dz\,&\leq \left(C\,+\,\frac{C}{\rho^{2}}\right)\left[\Vert Du\Vert_{L^{p}(Q_{R})}^{p}+\Vert Du\Vert_{L^{p}(Q_{R})}^{2}+\lambda^{p}+\lambda^{2}+1\right]\nonumber\\
&\,\,\,\,\,\,\,+\,C\,\Vert f\Vert_{L^{p'}\left(t_{0}-R^{2},t_{0};B_{p',1}^{\frac{p-2}{p}}(B_{R}(x_{0}))\right)}^{p'}
\end{align} for a positive constant $C$ depending only on $n$, $p$ and $R$
in the case $\lambda=0$, and additionally on $\alpha$ if $\lambda>0$.
Besides, if $\lambda=0$ we get 
\[
Du\,\in\,L_{loc}^{p}\left(0,T;W_{loc}^{\sigma,p}(\Omega,\mathbb{R}^{n})\right)\,\,\,\,\,\,\,\,\,\mathit{for}\,\,\mathit{all}\,\,\sigma\in\left(0,\frac{2}{p}\right).
\]
\end{thm}

\noindent \begin{brem} Notice that, for every $\alpha\geq0$, we
have 
\begin{equation}
\mathcal{V}_{\alpha,0}(\xi)=\,\frac{2}{p}\,V(\xi):=\,\frac{2}{p}\,\vert\xi\vert^{\frac{p-2}{2}}\xi\,,\label{eq:V0}
\end{equation}
i.e. $\mathcal{V}_{0}$ is actually independent of the parameter $\alpha$,
which explains the choice $(\ref{eq:alpha})_{1}$ in the statement
of the above theorem. The condition $(\ref{eq:alpha})_{2}$ will instead
be decisive to carry out the proof of Proposition \ref{prop:UniformSobolev}
below. Moreover, looking at (\ref{eq:V0}), one can easily understand
that, on the one hand, Theorem \ref{thm:theo1} extends the results
proved in \cite{Duzaar} to a widely degenerate parabolic setting.
On the other hand, it extends the aforementioned results to the case
of data in a suitable Lebesgue-Besov parabolic space, which turns
out to be optimal, as can be seen by appropriately adapting the example
in \cite[Section 5]{BraSan} to the parabolic setting (in this regard,
see also \cite[page 2]{AmGrPa}).\end{brem}\medskip{}

\noindent $\hspace*{1em}$As already mentioned, for $\lambda>0$ the
weak solutions of (\ref{eq:1}) may not be twice weakly differentiable
with respect to the spatial variable. Therefore, in general we cannot
differentiate the equation to estimate the second spatial derivatives
of the solutions. Here, we overcome this difficulty by introducing
a suitable family of approximating problems, whose solutions are regular
enough by standard parabolic regularity (see Section \ref{sec:a priori}
below). The major effort in proving Theorem \ref{thm:theo1} is to
establish suitable estimates for the solutions $u_{\varepsilon}$
of the regularized problems that are uniform with respect to the approximation
parameter $\varepsilon$. Next, we take advantage of these uniform
estimates in the use of a comparison argument aimed to transfer the
higher differentiability in space of $\mathcal{V}_{\alpha,\lambda}(Du_{\varepsilon})$
to the function $\mathcal{V}_{\alpha,\lambda}(Du)$.\\
The main differences with respect to the arguments used in \cite{AmPa}
lie in the use of the function $\mathcal{V}_{\alpha,\lambda}(Du)$
in place of $H_{p/2}(Du)$ and in the derivation of the \textit{a
priori} estimates for the terms involving the datum $f$, which now
only has a Besov spatial regularity. Indeed, if $f$ had a suitable
Sobolev spatial regularity as in \cite{AmPa}, the terms of the estimate
coming from the inhomogeneity of the equation could be controlled
using the information on the integrability of the spatial gradient
of the solution. Here, instead, in order to control the terms coming
from the inhomogeneity, we essentially use the duality of Besov spaces
and Theorem \ref{negder} below: this novelty comes from an idea that
we have already exploited in the recent paper \cite{AmGrPa}, in the
elliptic setting.\\
$\hspace*{1em}$Now we turn our attention to the case $p=2$. It is
well known that for the non-homogeneous heat equation 
\begin{equation}
\partial_{t}u-\mathrm{div}\,Du=f\,\,\,\,\,\,\,\,\mathrm{in}\,\,\,\Omega_{T},\label{eq:heat}
\end{equation}
the spatial $W^{1,2}$-regularity of the weak spatial gradient $Du$
can be achieved by simply assuming that $f\in L_{loc}^{2}(\Omega_{T})$.
Here, we prove that a similar result holds even when dealing with
much more degenerate equations. More precisely, we establish the following
result. 
\begin{thm}
\noindent \label{thm:nuovo} Let $n\geq2$, $\lambda\geq0$ and $f\in L_{loc}^{2}(\Omega_{T})$.
Moreover, let 
\[
\alpha=\begin{cases}
\begin{array}{cc}
0 & \textit{if}\,\,\,\lambda=0,\\
\textit{any\,\,value\,\,in}\,\,\left[\frac{3}{2},\infty\right) & \textit{if}\,\,\,\lambda>0,
\end{array}\end{cases}
\]
and assume that
\[
u\in C^{0}\left(0,T;L^{2}(\Omega)\right)\cap L_{loc}^{2}\left(0,T;W_{loc}^{1,2}(\Omega)\right)
\]
is a local weak solution of the equation 
\[
\partial_{t}u-\mathrm{div}\left((\vert Du\vert-\lambda)_{+}\,\frac{Du}{\vert Du\vert}\right)=f\,\,\,\,\,\,\,\,\mathit{in}\,\,\,\Omega_{T}.
\]
Then 
\[
\mathcal{V}_{\alpha,\lambda}(Du)\,\in\,L_{loc}^{2}\left(0,T;W_{loc}^{1,2}(\Omega,\mathbb{R}^{n})\right).
\]
Furthermore, for any parabolic cylinder $Q_{r}(z_{0})\subset Q_{\rho}(z_{0})\subset Q_{R}(z_{0})\Subset\Omega_{T}$
we have
\[
\int_{Q_{r/2}(z_{0})}\left|D_{x}\mathcal{V}_{\alpha,\lambda}(Du)\right|^{2}dz\,\leq\,\frac{C}{\rho^{2}}\left(\Vert Du\Vert_{L^{2}(Q_{R})}^{2}+\lambda^{2}+1\right)\,+\,C\left(\Vert f\Vert_{L^{2}(Q_{R})}^{2}+1\right)
\]
for a positive constant $C$ depending only on $n$ and $R$ in the
case $\lambda=0$, and additionally on $\alpha$ if $\lambda>0$.
\end{thm}

\noindent \begin{brem} In the case $p=2$, we have $\mathcal{V}_{0}(Du)=Du$.
Then, for $\lambda=0$, it follows from Theorem \ref{thm:nuovo} that
\[
Du\,\in\,L_{loc}^{2}\left(0,T;W_{loc}^{1,2}(\Omega,\mathbb{R}^{n})\right).
\]
As already mentioned, this is a well-established result for equation
(\ref{eq:heat}) when $f\in L_{loc}^{2}(\Omega_{T})$.\end{brem}

\noindent $\hspace*{1em}$For $p>2$, we now consider the evolutionary
$p$-Poisson equation
\begin{equation}
\partial_{t}u-\mathrm{div}\,(\vert Du\vert^{p-2}Du)=f\,\,\,\,\,\,\,\,\mathrm{in}\,\,\,\Omega_{T}\,,\label{eq:p-Poisson}
\end{equation}
which is obtained from equation (\ref{eq:1}) by setting $\lambda=0$.
As we anticipated earlier, from Theorem \ref{thm:theo1} we can easily
deduce that the weak solutions of (\ref{eq:p-Poisson}) admit a weak
time derivative which belongs to the local Lebesgue space $L_{loc}^{p'}(\Omega_{T})$.
The idea is roughly as follows. Since the above-mentioned theorem
tells us that in a certain pointwise sense the second spatial derivatives
of $u$ exist, we may develop the expression under the divergence
symbol. This will give us an expression that equals $\partial_{t}u$,
from which we get the desired integrability of the time derivative.
Such an argument must be made more rigorous. Furthermore, we also
need to make explicit\textit{ a priori} local estimates. These are
provided in the following result. 
\begin{thm}
\label{thm:theo2} Let $n\geq2$, $p>2$ and $f\in L_{loc}^{p'}\left(0,T;B_{p',1,loc}^{\frac{p-2}{p}}(\Omega)\right)$.
Moreover, assume that 
\[
u\in C^{0}\left(0,T;L^{2}(\Omega)\right)\cap L_{loc}^{p}\left(0,T;W_{loc}^{1,p}(\Omega)\right)
\]
is a local weak solution of equation $\mathrm{(\ref{eq:p-Poisson})}$.
Then, the time derivative of $u$ exists in the weak sense and satisfies
\[
\partial_{t}u\,\in\,L_{loc}^{p'}(\Omega_{T}).
\]
Furthermore, for any parabolic cylinder $Q_{r}(z_{0})\subset Q_{\rho}(z_{0})\subset Q_{R}(z_{0})\Subset\Omega_{T}$
we have \begin{align}\label{eq:stimaapriori2}
\left(\int_{Q_{r/2}(z_{0})}\left|\partial_{t}u\right|^{p'}dz\right)^{\frac{1}{p'}}& \leq \left(C\,+\,\frac{C}{\rho}\right)\left[\Vert Du\Vert_{L^{p}(Q_{R})}^{p-1}+\Vert Du\Vert_{L^{p}(Q_{R})}^{\frac{p}{2}}+\Vert Du\Vert_{L^{p}(Q_{R})}^{\frac{p-2}{2}}\right]\nonumber\\
&\,\,\,\,\,\,\,+\,C\,\Vert Du\Vert_{L^{p}(Q_{R})}^{\frac{p-2}{2}}\,\Vert f\Vert_{L^{p'}\left(t_{0}-R^{2},t_{0};B_{p',1}^{\frac{p-2}{p}}(B_{R}(x_{0}))\right)}^{\frac{p'}{2}}\,+\,\Vert f\Vert_{L^{p'}(Q_{R})}
\end{align}for a positive constant $C$ depending only on $n$, $p$ and $R$.
\end{thm}

\noindent $\hspace*{1em}$Before describing the structure of this
paper, we wish to point out that, starting from the weaker assumption
\[
f\in L_{loc}^{p'}\left(0,T;B_{p',1,loc}^{\frac{p-2}{p}}(\Omega)\right)\,\,\,\,\,\,\,\,\mathrm{with}\,\,\,p>2,
\]
Sobolev regularity results such as those of Theorems \ref{thm:theo1}
and \ref{thm:theo2} seem not to have been established yet for weak
solutions to parabolic PDEs that are far less degenerate than equation
(\ref{eq:1}) with $\lambda>0$. In particular, the results contained
in this paper permit to improve the existing literature, already for
the evolutionary $p$-Poisson equation (\ref{eq:p-Poisson}), which
exhibits a milder degeneracy. Moreover, as far as we know, for $\lambda>0$
the result of Theorem \ref{thm:nuovo} is completely new. 

\subsection{Plan of the paper }

$\hspace*{1em}$The paper is organized as follows. Section \ref{sec:prelim}
is devoted to the preliminaries: after a list of some classic notations
and some essentials estimates, we recall the fundamental properties
of the difference quotients of Sobolev functions. In Section \ref{sec:Functional spaces},
we recall the basic facts on the functional spaces involved in this
paper. Here, the most important points are Theorems \ref{duality}
and \ref{negder}, whose role is crucial in the proof of Theorem \ref{thm:theo1}.
In Section \ref{sec:a priori}, we consider a regularization of equation
(\ref{eq:1}) and establish some \textit{a priori} estimates for its
weak solution. In particular, we prove Sobolev estimates which are
independent of the regularization parameter (Propositions \ref{prop:UniformSobolev}
and \ref{prop:UniformSobolev-2}). These estimates will be needed
to demonstrate Theorems \ref{thm:theo1} and \ref{thm:nuovo}, whose
proofs are achieved in Section \ref{sec:proof1}. Finally, in Section
\ref{sec:proof2} we prove Theorem \ref{thm:theo2}.

\section{Preliminaries \label{sec:prelim}}

\subsection{Notation and essential definitions }

\noindent $\hspace*{1em}$In this paper we shall denote by $C$ or
$c$ a general positive constant that may vary on different occasions.
Relevant dependencies on parameters and special constants will be
suitably emphasized using parentheses or subscripts. The norm we use
on $\mathbb{R}^{k}$, $k\in\mathbb{N}$, will be the standard Euclidean
one and it will be denoted by $\left|\,\cdot\,\right|$. In particular,
for the vectors $\xi,\eta\in\mathbb{R}^{k}$, we write $\langle\xi,\eta\rangle$
for the usual inner product and $\left|\xi\right|:=\langle\xi,\xi\rangle^{\frac{1}{2}}$
for the corresponding Euclidean norm.\\
$\hspace*{1em}$For points in space-time, we will frequently use abbreviations
like $z=(x,t)$ or $z_{0}=(x_{0},t_{0})$, for spatial variables $x$,
$x_{0}\in\mathbb{R}^{n}$ and times $t$, $t_{0}\in\mathbb{R}$. We
also denote by $B_{\rho}(x_{0})=\left\{ x\in\mathbb{R}^{n}:\left|x-x_{0}\right|<\rho\right\} $
the $n$-dimensional open ball with radius $\rho>0$ and center $x_{0}\in\mathbb{R}^{n}$;
when not important, or clear from the context, we shall omit to denote
the center as follows: $B_{\rho}\equiv B_{\rho}(x_{0})$. Unless otherwise
stated, different balls in the same context will have the same center.
Moreover, we use the notation 
\[
Q_{\rho}(z_{0}):=B_{\rho}(x_{0})\times(t_{0}-\rho^{2},t_{0}),\,\,\,\,\,z_{0}=(x_{0},t_{0})\in\mathbb{R}^{n}\times\mathbb{R},\,\,\rho>0,
\]
for the backward parabolic cylinder with vertex $(x_{0},t_{0})$ and
width $\rho$. We shall sometimes omit the dependence on the vertex
when all cylinders occurring in a proof share the same vertex. For
a general cylinder $Q=B\times(t_{0},t_{1})$, where $B\subset\mathbb{R}^{n}$
and $t_{0}<t_{1}$, we denote by 
\[
\partial_{\mathrm{par}}Q:=(\overline{B}\times\{t_{0}\})\cup(\partial B\times(t_{0},t_{1}))
\]
the usual \textit{parabolic boundary} of $Q$, which is nothing but
its standard topological boundary without the upper cap $\overline{B}\times\{t_{1}\}$.\\
$\hspace*{1em}$If $E\subseteq\mathbb{R}^{k}$ is a Lebesgue-measurable
set, then we will denote by $\vert E\vert$ its $k$-dimensional Lebesgue
measure.\\
$\hspace*{1em}$For every function $v\in L_{loc}^{1}(Q,\mathbb{R}^{k})$,
where $Q\subset\mathbb{R}^{n+1}$ and $k\in\mathbb{N}$, we define
the mollified function $v_{\varrho}$ as follows:
\begin{equation}
v_{\varrho}(z):=\int_{\mathbb{R}^{n+1}}v(\tilde{z})\,\phi_{\varrho}(z-\tilde{z})\,d\tilde{z},\label{eq:mollified}
\end{equation}
where 
\[
\phi_{\varrho}(z):=\,\frac{1}{\varrho^{n+1}}\,\phi_{1}\left(\frac{z}{\varrho}\right),\,\,\,\,\,\,\,\,\,\,\varrho>0,
\]
with $\phi_{1}\in C_{0}^{\infty}(\mathcal{B}_{1}(0))$\footnote{Here $\mathcal{B}_{1}(0)$ denotes the $(n+1)$-dimensional open unit
ball centered at the origin.} denoting the standard, non-negative, radially symmetric mollifier
in $\mathbb{R}^{n+1}$. Obviously, here the function $v$ is meant
to be extended by the $k$-dimensional null vector outside $Q$.\\
$\hspace*{1em}$For further needs, we now define the auxiliary function
$H_{\gamma}:\mathbb{R}^{n}\rightarrow\mathbb{R}^{n}$ by 
\begin{equation}
H_{\gamma}(\xi):=\begin{cases}
\begin{array}{cc}
(\vert\xi\vert-\lambda)_{+}^{\gamma}\frac{\xi}{\left|\xi\right|} & \mathrm{if}\,\,\,\xi\neq0,\\
0 & \mathrm{if}\,\,\,\xi=0,
\end{array}\end{cases}\label{eq:Hfun}
\end{equation}

\noindent where $\gamma>0$ is a parameter. We conclude this first
part of the preliminaries by recalling the following definition.
\begin{defn}
\noindent Let $\lambda\geq0$. A function
\[
u\in C^{0}\left(0,T;L^{2}(\Omega)\right)\cap L_{loc}^{p}\left(0,T;W_{loc}^{1,p}(\Omega)\right)
\]
is a \textit{local weak solution} of equation (\ref{eq:1}) if and
only if, for any test function $\varphi\in C_{0}^{\infty}(\Omega_{T})$,
the following integral identity holds:
\begin{equation}
\int_{\Omega_{T}}\left(u\cdot\partial_{t}\varphi-\langle H_{p-1}(Du),D\varphi\rangle\right)\,dz\,=\,-\int_{\Omega_{T}}f\varphi\,dz.\label{eq:weaksol}
\end{equation}
\end{defn}

\subsection{Algebraic inequalities}

\noindent $\hspace*{1em}$In this section, we gather some relevant
algebraic inequalities that will be needed later on. The first result
follows from an elementary computation.
\begin{lem}
\noindent \label{lem:alge1} For $\xi,\eta\in\mathbb{R}^{n}\setminus\{0\}$,
we have 
\[
\left|\frac{\xi}{\vert\xi\vert}-\frac{\eta}{\vert\eta\vert}\right|\leq\mathrm{\frac{2}{\vert\eta\vert}\,\vert\xi-\eta\vert.}
\]
\end{lem}

\noindent We now recall the following estimate, whose proof can be
found in \cite[Chapter 12]{Lind}.
\begin{lem}
\noindent \label{lem:Lind} Let $p\in(2,\infty)$ and $k\in\mathbb{N}$.
Then, for every $\xi,\eta\in\mathbb{R}^{k}$ we get 
\[
\vert\xi-\eta\vert^{p}\leq\,C\left|\vert\xi\vert^{\frac{p-2}{2}}\xi-\vert\eta\vert^{\frac{p-2}{2}}\eta\right|^{2}
\]
for a constant $C\equiv C(p)>0$.
\end{lem}

\noindent $\hspace*{1em}$For the function $H_{p-1}$ defined by (\ref{eq:Hfun})
with $\gamma=p-1$, we record the following estimates, which can be
obtained by suitably modifying the proofs of \cite[Lemma 4.1]{Br}
and \cite[Lemma 2.8]{BDGP}, respectively.
\begin{lem}
\begin{singlespace}
\noindent \label{lem:Brasco} Let $2\leq p<\infty$. Then, for every
$\xi,\eta\in\mathbb{R}^{n}$ we get 
\[
\langle H_{p-1}(\xi)-H_{p-1}(\eta),\xi-\eta\rangle\,\geq\,\frac{4}{p^{2}}\left|H_{\frac{p}{2}}(\xi)-H_{\frac{p}{2}}(\eta)\right|^{2}.
\]
Moreover, if $\vert\eta\vert>\lambda$ we have 
\[
\langle H_{p-1}(\xi)-H_{p-1}(\eta),\xi-\eta\rangle\,\geq\,\frac{1}{2^{p+1}}\,\frac{(\vert\eta\vert-\lambda)^{p}}{\vert\eta\vert\,(\vert\xi\vert+\vert\eta\vert)}\,\vert\xi-\eta\vert^{2}.
\]
\end{singlespace}
\end{lem}

\noindent $\hspace*{1em}$The next result concerns the function $\mathcal{G}_{\alpha,\lambda}$
defined by (\ref{eq:Gfun}).
\begin{lem}
\noindent \label{lem:Glemma1} Let $p\in[2,\infty)$, $\lambda\geq0$
and $\alpha\geq0$. Then 
\[
\mathcal{G}_{\alpha,\lambda}(t)\le\,\frac{2}{p}\,t^{\frac{p}{2}}\left(\frac{t}{t+\lambda}\right)^{\frac{1+2\alpha}{2}}
\]
for every $t>0$.
\end{lem}

\noindent \begin{proof}[\bfseries{Proof}] Since the function 
\[
K(\omega):=\left(\frac{\omega}{\omega+\lambda}\right)^{\frac{1+2\alpha}{2}},\,\,\,\,\,\,\,\,\,\,\omega>0,
\]
is non-decreasing, for every $t>0$ we have 
\[
\mathcal{G}_{\alpha,\lambda}(t)=\int_{0}^{t}K(\omega)\,\omega^{\frac{p}{2}-1}\,d\omega\,\le\left(\frac{t}{t+\lambda}\right)^{\frac{1+2\alpha}{2}}\int_{0}^{t}\omega^{\frac{p}{2}-1}\,d\omega\,=\,\frac{2}{p}\,t^{\frac{p}{2}}\left(\frac{t}{t+\lambda}\right)^{\frac{1+2\alpha}{2}}.
\]
\end{proof}

\noindent Finally, the next lemma relates the function $\mathcal{V}_{\alpha,\lambda}(\xi)$
with $H_{p-1}(\xi)$.
\begin{lem}
\label{lem:vital} Let $p\in[2,\infty)$, $\lambda\geq0$ and $\alpha\geq0$.
Then, there exists a constant $C\equiv C(p)>0$ such that 
\begin{equation}
\vert\mathcal{V}_{\alpha,\lambda}(\xi)-\mathcal{V}_{\alpha,\lambda}(\eta)\vert^{2}\leq\,C\,\langle H_{p-1}(\xi)-H_{p-1}(\eta),\xi-\eta\rangle\label{eq:VH}
\end{equation}
for every $\xi,\eta\in\mathbb{R}^{n}$. 
\end{lem}

\noindent \begin{proof}[\bfseries{Proof}] For $\lambda=0$ this is
a well-known result (see, for example, \cite[Lemma 2.2]{AmPa}). Therefore,
from now on we shall assume that $\lambda>0$. We first note that
inequality (\ref{eq:VH}) is trivially satisfied when $\vert\xi\vert,\vert\eta\vert\leq\lambda$.
If $\vert\eta\vert\leq\lambda<\vert\xi\vert$, using the definitions
(\ref{eq:Vfun}), (\ref{eq:Gfun}), (\ref{eq:Hfun}) and Lemma \ref{lem:Brasco},
we obtain\begin{align*}
\vert\mathcal{V}_{\alpha,\lambda}(\xi)-\mathcal{V}_{\alpha,\lambda}(\eta)\vert^{2}&=[\mathcal{G}_{\alpha,\lambda}(\vert\xi\vert-\lambda)]^{2}\le\left(\int_{0}^{\vert\xi\vert-\lambda}\omega^{\frac{p}{2}-1}\,d\omega\right)^{2}=\,\frac{4}{p^{2}}\,(\vert\xi\vert-\lambda)^{p}\\
&=\,\frac{4}{p^{2}}\left|H_{\frac{p}{2}}(\xi)-H_{\frac{p}{2}}(\eta)\right|^{2}\le\,\langle H_{p-1}(\xi)-H_{p-1}(\eta),\xi-\eta\rangle.
\end{align*}Now let $\vert\xi\vert,\vert\eta\vert>\lambda$. Without loss of generality,
we may assume that $\vert\eta\vert\geq\vert\xi\vert>\lambda$. This
implies 
\begin{equation}
\vert\eta\vert^{2}=\,\frac{\vert\eta\vert\,(\vert\eta\vert+\vert\eta\vert)}{2}\,\geq\,\frac{\vert\eta\vert\,(\vert\xi\vert+\vert\eta\vert)}{2}\,.\label{eq:implication}
\end{equation}
Moreover, we have \begin{align*}
\vert\mathcal{V}_{\alpha,\lambda}(\xi)-\mathcal{V}_{\alpha,\lambda}(\eta)\vert\,&=\,\left|\mathcal{V}_{\alpha,\lambda}(\xi)-\mathcal{G}_{\alpha,\lambda}(\vert\eta\vert-\lambda)\,\frac{\xi}{\left|\xi\right|}\,+\mathcal{G}_{\alpha,\lambda}(\vert\eta\vert-\lambda)\,\frac{\xi}{\left|\xi\right|}\,-\mathcal{V}_{\alpha,\lambda}(\eta)\right|\\
&\leq\,\left|\mathcal{G}_{\alpha,\lambda}(\vert\xi\vert-\lambda)-\mathcal{G}_{\alpha,\lambda}(\vert\eta\vert-\lambda)\right|\,+\,\mathcal{G}_{\alpha,\lambda}(\vert\eta\vert-\lambda)\left|\frac{\xi}{\left|\xi\right|}-\frac{\eta}{\left|\eta\right|}\right|\\
&\le\,\int_{\vert\xi\vert-\lambda}^{\vert\eta\vert-\lambda}\frac{\omega^{\frac{p-1+2\alpha}{2}}}{(\omega+\lambda)^{\frac{1+2\alpha}{2}}}\,d\omega\,+\,\frac{2}{\vert\eta\vert}\,\vert\xi-\eta\vert\int_{0}^{\vert\eta\vert-\lambda}\frac{\omega^{\frac{p-1+2\alpha}{2}}}{(\omega+\lambda)^{\frac{1+2\alpha}{2}}}\,d\omega\\
&\le\,\int_{\vert\xi\vert-\lambda}^{\vert\eta\vert-\lambda}\omega^{\frac{p}{2}-1}\,d\omega\,+\,\frac{2}{\vert\eta\vert}\,\vert\xi-\eta\vert\int_{0}^{\vert\eta\vert-\lambda}\omega^{\frac{p}{2}-1}\,d\omega\\
&=\,\frac{2}{p}\,\left||H_{\frac{p}{2}}(\eta)|-|H_{\frac{p}{2}}(\xi)|\right|+\,4\,\frac{(\vert\eta\vert-\lambda)^{\frac{p}{2}}}{p\,\vert\eta\vert}\,\vert\xi-\eta\vert\\
&\le\,\frac{2}{p}\,\left|H_{\frac{p}{2}}(\xi)-H_{\frac{p}{2}}(\eta)\right|+\,4\,\frac{(\vert\eta\vert-\lambda)^{\frac{p}{2}}}{p\,\vert\eta\vert}\,\vert\xi-\eta\vert,
\end{align*}where, in the third line, we have used Lemma \ref{lem:alge1} and
the fact that $\mathcal{G}_{\alpha,\lambda}$ is an increasing function.
Now, applying Young's inequality, estimate (\ref{eq:implication})
and Lemma \ref{lem:Brasco}, we obtain\begin{align*}
\vert\mathcal{V}_{\alpha,\lambda}(\xi)-\mathcal{V}_{\alpha,\lambda}(\eta)\vert^{2}\,&\le\,\frac{8}{p^{2}}\,\left|H_{\frac{p}{2}}(\xi)-H_{\frac{p}{2}}(\eta)\right|^{2}+\,32\,\frac{(\vert\eta\vert-\lambda)^{p}}{p^{2}\,\vert\eta\vert^{2}}\,\vert\xi-\eta\vert^{2}\\
&\le\,\frac{8}{p^{2}}\,\left|H_{\frac{p}{2}}(\xi)-H_{\frac{p}{2}}(\eta)\right|^{2}+\,\frac{64}{p^{2}}\,\frac{(\vert\eta\vert-\lambda)^{p}}{\vert\eta\vert\,(\vert\xi\vert+\vert\eta\vert)}\,\vert\xi-\eta\vert^{2}\\
&\leq\,C_{p}\,\langle H_{p-1}(\xi)-H_{p-1}(\eta),\xi-\eta\rangle.
\end{align*}This completes the proof.\end{proof}

\subsection{Difference quotients \label{subsec:DiffOpe}}

\noindent $\hspace*{1em}$We recall here the definition and some elementary
properties of the difference quotients that will be useful in the
following (see, for example, \cite{Giu}).
\begin{defn}
\noindent For every vector-valued function $F:\mathbb{R}^{n}\rightarrow\mathbb{R}^{k}$
the \textit{finite difference operator }in the direction $x_{j}$
is defined by
\[
\tau_{j,h}F(x)=F(x+he_{j})-F(x),
\]
where $h\in\mathbb{R}$, $e_{j}$ is the unit vector in the direction
$x_{j}$ and $j\in\{1,\ldots,n\}$. \\
$\hspace*{1em}$The \textit{difference quotient} of $F$ with respect
to $x_{j}$ is defined for $h\in\mathbb{R}\setminus\{0\}$ by 
\[
\Delta_{j,h}F(x)\,=\,\frac{\tau_{j,h}F(x)}{h}\,.
\]
\end{defn}

\noindent When no confusion can arise, we shall omit the index $j$
and simply write $\tau_{h}$ or $\Delta_{h}$ instead of $\tau_{j,h}$
or $\Delta_{j,h}$, respectively. 
\begin{prop}
\noindent Let $\Omega\subset\mathbb{R}^{n}$ be an open set and let
$F\in W^{1,q}(\Omega)$, with $q\geq1$. Moreover, let $G:\Omega\rightarrow\mathbb{R}$
be a measurable function and consider the set
\[
\Omega_{\vert h\vert}:=\left\{ x\in\Omega:\mathrm{dist}\left(x,\partial\Omega\right)>\vert h\vert\right\} .
\]
Then:\\
\\
$\mathrm{(}\mathrm{i}\mathrm{)}$ $\Delta_{h}F\in W^{1,q}\left(\Omega_{\vert h\vert}\right)$
and $D_{x_{i}}(\Delta_{h}F)=\Delta_{h}(D_{x_{i}}F)$ for every $\,i\in\{1,\ldots,n\}$.\\

\noindent $\mathrm{(}\mathrm{ii}\mathrm{)}$ If at least one of the
functions $F$ or $G$ has support contained in $\Omega_{\vert h\vert}$,
then
\[
\int_{\Omega}F\,\Delta_{h}G\,dx\,=\,-\int_{\Omega}G\,\Delta_{-h}F\,dx.
\]
$\mathrm{(}\mathrm{iii}\mathrm{)}$ We have 
\[
\Delta_{h}(FG)(x)=F(x+he_{j})\Delta_{h}G(x)\,+\,G(x)\Delta_{h}F(x).
\]
\end{prop}

\noindent The next result about the finite difference operator is
a kind of integral version of the Lagrange Theorem and its proof can
be found in \cite[Lemma 8.1]{Giu}. 
\begin{lem}
\noindent \label{lem:Giusti1} If $0<\rho<R$, $\vert h\vert<\frac{R-\rho}{2}$,
$1<q<+\infty$ and $F\in L^{q}(B_{R},\mathbb{R}^{k})$ is such that
$DF\in L^{q}(B_{R},\mathbb{R}^{k\times n})$, then
\[
\int_{B_{\rho}}\left|\tau_{h}F(x)\right|^{q}dx\,\leq\,c^{q}(n)\,\vert h\vert^{q}\int_{B_{R}}\left|DF(x)\right|^{q}dx.
\]
Moreover 
\[
\int_{B_{\rho}}\left|F(x+he_{j})\right|^{q}dx\,\leq\,\int_{B_{R}}\left|F(x)\right|^{q}dx.
\]
\end{lem}

\noindent Finally, we recall the following fundamental result, whose
proof can be found in \cite[Lemma 8.2]{Giu}:
\begin{lem}
\noindent \label{lem:RappIncre} Let $F:\mathbb{R}^{n}\rightarrow\mathbb{R}^{k}$,
$F\in L^{q}(B_{R},\mathbb{R}^{k})$ with $1<q<+\infty$. Suppose that
there exist $\rho\in(0,R)$ and a constant $M>0$ such that 
\[
\sum_{j=1}^{n}\int_{B_{\rho}}\left|\tau_{j,h}F(x)\right|^{q}dx\,\leq\,M^{q}\,\vert h\vert^{q}
\]
for every $h$ with $\vert h\vert<\frac{R-\rho}{2}$. Then $F\in W^{1,q}(B_{\rho},\mathbb{R}^{k})$.
Moreover 
\[
\Vert DF\Vert_{L^{q}(B_{\rho})}\leq M
\]
and
\[
\Delta_{j,h}F\rightarrow D_{x_{j}}F\,\,\,\,\,\,\,\,\,\,\mathit{in}\,\,L_{loc}^{q}(B_{R},\mathbb{R}^{k})\,\,\,\,\mathit{as}\,\,h\rightarrow0,
\]
for each $j\in\{1,\ldots,n\}$.
\end{lem}

\selectlanguage{english}%

\section{Functional spaces \label{sec:Functional spaces}}

\noindent $\hspace*{1em}$Here we recall some essential facts about
the functional spaces involved in this paper, starting with the definition
and some properties of the Besov spaces that will be useful to prove
our results.\\
We denote by $\mathcal{S}(\mathbb{R}^{n})$ and $\mathcal{S}'(\mathbb{R}^{n})$
the Schwartz space and the space of tempered distributions on $\mathbb{R}^{n}$,
respectively. If $v\in\mathcal{S}(\mathbb{R}^{n})$, then 
\begin{equation}
\hat{v}(\xi)=(\mathcal{F}v)(\xi)=(2\pi)^{-n/2}\int_{\mathbb{R}^{n}}e^{-i\,\langle x,\xi\rangle}\,v(x)\,dx,\,\,\,\,\,\,\,\,\xi\in\mathbb{R}^{n},\label{eq:Fourier}
\end{equation}
denotes the Fourier transform of $v$. As usual, $\mathcal{F}^{-1}v$
and $v^{\vee}$ stand for the inverse Fourier transform, given by
the right-hand side of (\ref{eq:Fourier}) with $i$ in place of $-i$.
Both $\mathcal{F}$ and $\mathcal{F}^{-1}$ are extended to $\mathcal{S}'(\mathbb{R}^{n})$
in the standard way.\\
Now, let $\Gamma(\mathbb{R}^{n})$ be the collection of all sequences
$\varphi=\{\varphi_{j}\}_{j=0}^{\infty}\subset\mathcal{S}(\mathbb{R}^{n})$
such that 
\[
\begin{cases}
\begin{array}{cc}
\mathrm{supp}\,\varphi_{0}\subset\{x\in\mathbb{R}^{n}:|x|\le2\}\quad\quad\quad\quad\,\,\,\\
\mathrm{supp}\,\varphi_{j}\subset\{x\in\mathbb{R}^{n}:2^{j-1}\le|x|\le2^{j+1}\} & \mathrm{if}\,\,j\in\mathbb{N},
\end{array}\end{cases}
\]
for every multi-index $\beta$ there exists a positive number $c_{\beta}$
such that 
\[
2^{j\vert\beta\vert}\,\vert D^{\beta}\varphi_{j}(x)\vert\le c_{\beta}\,,\,\,\,\,\,\,\,\,\forall\,j\in\mathbb{N}_{0},\ \forall\,x\in\mathbb{R}^{n}
\]
and 
\[
\sum_{j=0}^{\infty}\varphi_{j}(x)=1\,,\,\,\,\,\,\,\,\,\forall\,x\in\mathbb{R}^{n}.
\]
Then, it is well known that $\Gamma(\mathbb{R}^{n})$ is not empty
(see \cite[Section 2.3.1, Remark 1]{Triebel}). Moreover, if $\{\varphi_{j}\}_{j=0}^{\infty}\in\Gamma(\mathbb{R}^{n})$,
the entire analytic functions $(\varphi_{j}\,\hat{v})^{\vee}(x)$
make sense pointwise in $\mathbb{R}^{n}$ for any $v\in\mathcal{S}'(\mathbb{R}^{n})$.
Therefore, the following definition makes sense: 
\selectlanguage{british}%
\begin{defn}
\noindent Let $s\in\mathbb{R}$, $1\leq p,q\le\infty$ and $\varphi=\{\varphi_{j}\}_{j=0}^{\infty}\in\Gamma(\mathbb{R}^{n})$.
We define the \textit{Besov space} $B_{p,q}^{s}(\mathbb{R}^{n})$
as the set of all $v\in\mathcal{S}'(\mathbb{R}^{n})$ such that 
\begin{equation}
\Vert v\Vert_{B_{p,q}^{s}(\mathbb{R}^{n})}:=\left(\sum_{j=0}^{\infty}2^{jsq}\,\Vert(\varphi_{j}\,\hat{v})^{\vee}\Vert_{L^{p}(\mathbb{R}^{n})}^{q}\right)^{\frac{1}{q}}<+\infty\,\,\,\,\,\,\,\,\,\,\mathrm{if}\,\,q<\infty,\label{eq:quasi-norm}
\end{equation}
and 
\begin{equation}
\Vert v\Vert_{B_{p,q}^{s}(\mathbb{R}^{n})}:=\,\sup_{j\,\in\,\mathbb{N}_{0}}\,2^{js}\,\Vert(\varphi_{j}\,\hat{v})^{\vee}\Vert_{L^{p}(\mathbb{R}^{n})}<+\infty\,\,\,\,\,\,\,\,\,\,\mathrm{if}\,\,q=\infty.\label{eq:quasi-norm2}
\end{equation}

\selectlanguage{english}%
\noindent \begin{brem}\foreignlanguage{british}{ The space $B_{p,q}^{s}(\mathbb{R}^{n})$
defined above is a Banach space with respect to the norm $\Vert\cdot\Vert_{B_{p,q}^{s}(\mathbb{R}^{n})}$.
Obviously, this norm depends on the chosen sequence $\varphi\in\Gamma(\mathbb{R}^{n})$,
but this is not the case for the spaces $B_{p,q}^{s}(\mathbb{R}^{n})$
themselves, in the sense that two different choices for the sequence
$\varphi$ give rise to equivalent norms (see \cite[Sections 2.3.2 and 2.3.3]{Triebel}).
This justifies our omission of the dependence on $\varphi$ in the
left-hand side of (\ref{eq:quasi-norm})$-$(\ref{eq:quasi-norm2})
and in the sequel.}\end{brem}
\end{defn}

\noindent $\hspace*{1em}$The norms of the \textit{classical Besov
spaces} $B_{p,q}^{s}(\mathbb{R}^{n})$ with $s\in(0,1)$, $1\leq p<\infty$
and $1\le q\le\infty$ can be characterized via differences of the
functions involved, cf. \cite[Section 2.5.12, Theorem 1]{Triebel}.
More precisely, for $h\in\mathbb{R}^{n}$ and a measurable function
$v:\mathbb{R}^{n}\rightarrow\mathbb{R}^{k}$, let us define 
\[
\delta_{h}v(x):=\,v(x+h)-v(x).
\]
Then we have the equivalence 
\begin{equation}
\Vert v\Vert_{B_{p,q}^{s}(\mathbb{R}^{n})}\,\approx\,\Vert v\Vert_{L^{p}(\mathbb{R}^{n})}\,+\,[v]_{B_{p,q}^{s}(\mathbb{R}^{n})}\,,\label{eq:equivalence}
\end{equation}

\noindent where 
\begin{equation}
[v]_{B_{p,q}^{s}(\mathbb{R}^{n})}:=\biggl({\displaystyle \int_{\mathbb{R}^{n}}\biggl({\displaystyle \int_{\mathbb{R}^{n}}\dfrac{|\delta_{h}v(x)|^{p}}{|h|^{sp}}\,dx\biggr)^{\frac{q}{p}}\dfrac{dh}{|h|^{n}}\biggr)^{\frac{1}{q}},\,\,\,\,\,\,\,\,\,\,\text{if}\,\,\,1\le q<\infty},}\label{eq:BeNorm1}
\end{equation}

\noindent and 
\begin{equation}
[v]_{B_{p,\infty}^{s}(\mathbb{R}^{n})}:={\displaystyle \sup_{h\,\in\,\mathbb{R}^{n}}\biggl({\displaystyle \int_{\mathbb{R}^{n}}\dfrac{|\delta_{h}v(x)|^{p}}{|h|^{sp}}\,dx\biggr)^{\frac{1}{p}}}}.\label{eq:BeNorm2}
\end{equation}

\noindent In (\ref{eq:BeNorm1}), if one simply integrates for $\vert h\vert<r$
for a fixed $r>0$, then an equivalent norm is obtained, since 
\begin{center}
$\biggl({\displaystyle \int_{\{|h|\,\geq\,r\}}\biggl({\displaystyle \int_{\mathbb{R}^{n}}\dfrac{|\delta_{h}v(x)|^{p}}{|h|^{sp}}\,dx\biggr)^{\frac{q}{p}}\dfrac{dh}{|h|^{n}}\biggr)^{\frac{1}{q}}\leq\,c(n,s,p,q,r)\,\Vert v\Vert_{L^{p}(\mathbb{R}^{n})}}}\,.$ 
\par\end{center}

\noindent Similarly, in (\ref{eq:BeNorm2}) one can simply take the
supremum over $|h|\leq r$ and obtain an equivalent norm. By construction,
$B_{p,q}^{s}(\mathbb{R}^{n})\subset L^{p}(\mathbb{R}^{n})$.\foreignlanguage{english}{}\\
\foreignlanguage{english}{$\hspace*{1em}$}Let $\Omega$ be an arbitrary
open set in $\mathbb{R}^{n}$. As usual, $\mathcal{D}(\Omega)=C_{0}^{\infty}(\Omega)$
stands for the space of all infinitely differentiable functions in
$\mathbb{R}^{n}$ with compact support in $\Omega$. Let $\mathcal{D}'(\Omega)$
be the dual space of all distributions in $\Omega$ and let $g\in\mathcal{S}'(\mathbb{R}^{n})$.
Then we denote by $g\vert_{\Omega}$ its restriction to $\Omega$,
i.e. 
\[
g\vert_{\Omega}\in\mathcal{D}'(\Omega):\,\,\,\,\,(g\vert_{\Omega})(\phi)=g(\phi)\,\,\,\,\,\,\,\mathrm{for}\,\,\phi\in\mathcal{D}(\Omega).
\]

\begin{defn}
\noindent Let $\Omega$ be an arbitrary domain in $\mathbb{R}^{n}$
with $\Omega\neq\mathbb{R}^{n}$ and let $s\in\mathbb{R}$, $1\leq p\le\infty$
and $1\leq q\le\infty$. Then 
\[
B_{p,q}^{s}(\Omega):=\left\{ v\in\mathcal{D}'(\Omega):\,v=g\vert_{\Omega}\,\,\,\mathrm{for\,\,some}\,\,g\in B_{p,q}^{s}(\mathbb{R}^{n})\right\} 
\]
and 
\[
\Vert v\Vert_{B_{p,q}^{s}(\Omega)}:=\,\inf\,\Vert g\Vert_{B_{p,q}^{s}(\mathbb{R}^{n})}\,,
\]
where the infimum is taken over all $g\in B_{p,q}^{s}(\mathbb{R}^{n})$
such that $g\vert_{\Omega}=v$. 
\end{defn}

\noindent $\hspace*{1em}$If $\Omega$ is a bounded $C^{\infty}$-domain
in $\mathbb{R}^{n}$, then the \textit{restriction operator} 
\[
\mathrm{re}_{\Omega}:\mathcal{S}'(\mathbb{R}^{n})\hookrightarrow\mathcal{D}'(\Omega),\,\,\,\,\,\,\mathrm{re}_{\Omega}(v)=v\vert_{\Omega}
\]
generates a linear and bounded map from $B_{p,q}^{s}(\mathbb{R}^{n})$
onto $B_{p,q}^{s}(\Omega)$. Furthermore, the spaces $B_{p,q}^{s}(\Omega)$
satisfy the so-called \textit{extension property}, as ensured by the
next theorem. 
\selectlanguage{english}%
\begin{thm}
\noindent \label{thm:extension} \foreignlanguage{british}{Let $s\in\mathbb{R}$,
let $1\leq p,q\le\infty$ and let $\Omega$ be a bounded $C^{\infty}$-domain
in $\mathbb{R}^{n}$. Then, there exists a linear and bounded extension
operator $\mathrm{ext}_{\Omega}:B_{p,q}^{s}(\Omega)\hookrightarrow B_{p,q}^{s}(\mathbb{R}^{n})$
such that $\mathrm{re}_{\Omega}\circ\mathrm{ext}_{\Omega}=\mathrm{id}$,
where $\mathrm{id}$ is the identity in $B_{p,q}^{s}(\Omega)$}.
\end{thm}

\selectlanguage{british}%
\noindent We refer to \cite[Theorem 2.82]{TriebelIV} for a proof
of the previous theorem (see also \cite[Theorem 3.3.4]{Triebel}).\\
$\hspace*{1em}$If $s\in\mathbb{R}$, $1\leq p<\infty$ and $1\leq q<\infty$,
then $\mathcal{S}(\mathbb{R}^{n})$ is a dense subset of $B_{p,q}^{s}(\mathbb{R}^{n})$
(cf. \cite[Theorem 2.3.3]{Triebel}). Consequently, in that case,
a continuous linear functional on $B_{p,q}^{s}(\mathbb{R}^{n})$ can
be interpreted in the usual way as an element of $\mathcal{S}'(\mathbb{R}^{n})$.
More precisely, $g\in\mathcal{S}'(\mathbb{R}^{n})$ belongs to the
dual space $(B_{p,q}^{s}(\mathbb{R}^{n}))'$ of the space $B_{p,q}^{s}(\mathbb{R}^{n})$
if and only if there exists a positive number $c$ such that 
\[
\vert g(\phi)\vert\leq\,c\,\Vert\phi\Vert_{B_{p,q}^{s}(\mathbb{R}^{n})}\,\,\,\,\,\,\,\,\,\,\mathrm{for\,\,all\,\,}\phi\in\mathcal{S}(\mathbb{R}^{n})\,.
\]
We endow $(B_{p,q}^{s}(\mathbb{R}^{n}))'$ with the natural dual norm,
defined by 
\[
\Vert g\Vert_{(B_{p,q}^{s}(\mathbb{R}^{n}))'}=\,\sup\,\left\{ \vert g(\phi)\vert:\phi\in\mathcal{S}(\mathbb{R}^{n})\,\,\,\mathrm{and}\,\,\,\Vert\phi\Vert_{B_{p,q}^{s}(\mathbb{R}^{n})}\leq1\right\} ,\,\,\,\,\,\,\,\,g\in(B_{p,q}^{s}(\mathbb{R}^{n}))'.
\]
Now we recall the following duality formula, which has to be meant
as an isomorphism of normed spaces (see \cite[Theorem 2.11.2]{Triebel}). 
\begin{thm}
\noindent \label{thm:duality00} Let $s\in\mathbb{R}$, $1\le p<\infty$
and $1\le q<\infty$. Then 
\[
(B_{p,q}^{s}(\mathbb{R}^{n}))'=B_{p',q'}^{-s}(\mathbb{R}^{n})\,,
\]
where $p'=\infty$ if $p=1$ (similarly for $q'$). 
\end{thm}

\noindent \begin{brem} The restrictions $p<\infty$ and $q<\infty$
in Theorem \ref{thm:duality00} are natural, since, if either $p=\infty$
or $q=\infty$, then $\mathcal{S}(\mathbb{R}^{n})$ is not dense in
$B_{p,q}^{s}(\mathbb{R}^{n})$, and the density of $\mathcal{S}(\mathbb{R}^{n})$
in $B_{p,q}^{s}(\mathbb{R}^{n})$ is the basis of our interpretation
of the dual space $(B_{p,q}^{s}(\mathbb{R}^{n}))'$.\end{brem}\smallskip{}

\noindent $\hspace*{1em}$For our purposes, we now give the following
definition. 
\begin{defn}
\noindent For $s\in\mathbb{R}$, $1\le p\leq\infty$ and $1\le q\leq\infty$,
we define $\mathring{B}_{p,q}^{s}(\mathbb{R}^{n})$ as the completion
of $\mathcal{S}(\mathbb{R}^{n})$ in $B_{p,q}^{s}(\mathbb{R}^{n})$
with respect to the norm 
\[
v\mapsto\Vert v\Vert_{B_{p,q}^{s}(\mathbb{R}^{n})}\,.
\]
Of course, only the limit cases $\max\,\{p,q\}=\infty$ are of interest.
We shall denote by $(\mathring{B}_{p,q}^{s}(\mathbb{R}^{n}))'$ the
topological dual of\textit{ }$\mathring{B}_{p,q}^{s}(\mathbb{R}^{n})$,
which is endowed with the natural dual norm 
\[
\Vert g\Vert_{(\mathring{B}_{p,q}^{s}(\mathbb{R}^{n}))'}=\,\sup\,\left\{ \vert g(\phi)\vert:\phi\in\mathcal{S}(\mathbb{R}^{n})\,\,\,\mathrm{and}\,\,\,\Vert\phi\Vert_{B_{p,q}^{s}(\mathbb{R}^{n})}\leq1\right\} ,\,\,\,\,\,\,\,\,g\in(\mathring{B}_{p,q}^{s}(\mathbb{R}^{n}))'.
\]
\end{defn}

\noindent $\hspace*{1em}$The following duality result can be found
in \cite[Section 2.11.2, Remark 2]{Triebel} (see also \cite[pages 121 and 122]{Tri0}). 
\begin{thm}
\label{duality} Let $s\in\mathbb{R}$, $1\le p\le\infty$ and $1\le q\le\infty$.
Then 
\[
(\mathring{B}_{p,q}^{s}(\mathbb{R}^{n}))'=B_{p',q'}^{-s}(\mathbb{R}^{n})\,,
\]
where $p'=1$ if $p=\infty$ (similarly for $q'$). 
\end{thm}

\noindent $\hspace*{1em}$The next result is a key ingredient for
the proof of Theorem \ref{thm:theo1} and its proof can be found in
\cite[Section 3.3.5]{Triebel}. 
\selectlanguage{english}%
\begin{thm}
\label{negder} \foreignlanguage{british}{Let $s\in\mathbb{R}$ and
$1\leq p,q\le\infty$. Moreover, assume that $\Omega$ is a bounded
$C^{\infty}$-domain in $\mathbb{R}^{n}$. Then, for every $v\in B_{p,q}^{s}(\Omega)$
and every $j\in\{1,\ldots,n\}$ we have 
\[
\Vert\partial_{j}v\Vert_{B_{p,q}^{s-1}(\Omega)}\,\le c\,\Vert v\Vert_{B_{p,q}^{s}(\Omega)}
\]
for a positive constant $c$ which is independent of $v$}.
\end{thm}

\noindent $\hspace*{1em}$We can also define local Besov spaces as
follows. Given a domain $\Omega\subset\mathbb{R}^{n}$, we say that
a function $v$ belongs to $B_{p,q,loc}^{s}(\Omega)$ if $\phi\,v\in B_{p,q}^{s}(\mathbb{R}^{n})$
whenever $\phi\in C_{0}^{\infty}(\Omega)$.
\selectlanguage{british}%
\begin{defn}
\noindent Let $\Omega\subseteq\mathbb{R}^{n}$ be an open set. For
any $s\in(0,1)$ and for any $q\in[1,+\infty)$, we define the \textit{fractional
Sobolev space} $W^{s,q}(\Omega,\mathbb{R}^{k})$ as follows: 
\[
W^{s,q}(\Omega,\mathbb{R}^{k}):=\left\{ v\in L^{q}(\Omega,\mathbb{R}^{k}):\frac{\left|v(x)-v(y)\right|}{\left|x-y\right|^{\frac{n}{q}\,+\,s}}\,\in L^{q}\left(\Omega\times\Omega\right)\right\} ,
\]
i.e. an intermerdiate Banach space between $L^{q}(\Omega,\mathbb{R}^{k})$
and $W^{1,q}(\Omega,\mathbb{R}^{k})$, endowed with the norm 
\[
\Vert v\Vert_{W^{s,q}(\Omega)}:=\,\Vert v\Vert_{L^{q}(\Omega)}\,+\,[v]_{W^{s,q}(\Omega)}\,,
\]
where the term 
\begin{equation}
[v]_{W^{s,q}(\Omega)}:=\left(\int_{\Omega}\int_{\Omega}\frac{\left|v(x)-v(y)\right|^{q}}{\left|x-y\right|^{n\,+\,sq}}\,dx\,dy\right)^{\frac{1}{q}}\label{eq:Gagliardo}
\end{equation}
is the so-called \textit{Gagliardo seminorm} of $v$\foreignlanguage{english}{.}
\end{defn}

\noindent \begin{brem} For every $s\in(0,1)$ and every $q\in[1,\infty)$,
we have $B_{q,q}^{s}(\mathbb{R}^{n})=W^{s,q}(\mathbb{R}^{n})$. In
fact, using the change of variable $y=x+h$ in (\ref{eq:Gagliardo})
with $\Omega=\mathbb{R}^{n}$, one gets the seminorm (\ref{eq:BeNorm1})
with $p=q$.\end{brem}

\selectlanguage{english}%
\noindent $\hspace*{1em}$The following embedding result can be obtained
by combining \cite[Section 2.2.2, Remark 3]{Triebel} with \cite[Section 2.3.2, Proposition 2(ii)]{Triebel}.
\begin{lem}
\label{lem:emb} Let $s\in(0,1)$ and $q\ge1$. Then, for every $\sigma\in(0,1-s)$
we have the continuous embedding $W_{loc}^{s+\sigma,q}(\mathbb{R}^{n})\hookrightarrow B_{q,1,loc}^{s}(\mathbb{R}^{n})$. 
\end{lem}

\selectlanguage{british}%
\noindent $\hspace*{1em}$For the treatment of parabolic equations,
we now give the following definitions. 
\begin{defn}
\noindent Let $1\leq q<\infty$ and $0<\beta<1$. A map $g\in L^{q}(\Omega\times(t_{0},t_{1}),\mathbb{R}^{k})$
belongs to the space $L^{q}(t_{0},t_{1};W^{\beta,q}(\Omega,\mathbb{R}^{k}))$
if and only if
\[
\int_{t_{0}}^{t_{1}}\int_{\Omega}\int_{\Omega}\frac{\left|g(x,t)-g(y,t)\right|^{q}}{\left|x-y\right|^{n\,+\,\beta q}}\,dx\,dy\,dt<\infty.
\]
\end{defn}

\begin{defn}
Let $1\leq p,q<\infty$ and $0<s<1$. A map $g\in L^{p}(\Omega\times(t_{0},t_{1}),\mathbb{R}^{k})$
belongs to the space $L^{p}(t_{0},t_{1};B_{p,q}^{s}(\Omega,\mathbb{R}^{k}))$
if and only if
\[
\int_{t_{0}}^{t_{1}}\Vert g(\cdot,t)\Vert_{B_{p,q}^{s}(\Omega)}^{p}\,dt<\infty.
\]
\end{defn}

\noindent In this paper, we use the corresponding local versions of
the above spaces, which are denoted by the subscript ``\textit{loc}''.
More precisely, we write $g\in L_{loc}^{q}(0,T;W_{loc}^{\beta,q}(\Omega,\mathbb{R}^{k}))$
if and only if $g\in L^{q}(t_{0},t_{1};W^{\beta,q}(\Omega',\mathbb{R}^{k}))$
for all domains $\Omega'\times(t_{0},t_{1})\Subset\Omega_{T}$. The
local Lebesgue-Besov space $L_{loc}^{p}(0,T;B_{p,q,loc}^{s}(\Omega,\mathbb{R}^{k}))$
is defined in a similar way. Furthermore, we shall also use the following
notation, which is typical of Bochner spaces: 

\[
\Vert g\Vert_{L^{p}(t_{0},t_{1};B_{p,q}^{s}(\Omega'))}:=\left(\int_{t_{0}}^{t_{1}}\Vert g(\cdot,t)\Vert_{B_{p,q}^{s}(\Omega')}^{p}\,dt\right)^{\frac{1}{p}}.
\]

\noindent $\hspace*{1em}$We conclude this section with the parabolic
version of the well-known relation between Nikolskii spaces and fractional
Sobolev spaces (see \cite[7.73]{Adams}). This version is contained
in \cite[Lemma 2.4]{Sche} (see also \cite[Proposition 2.19]{Duzaar}),
and its proof can be obtained by a straightforward adaptation of the
standard elliptic results contained in \cite{DuGaMi,EspLeoMin,Ming1,Ming2,Ming3}. 
\begin{prop}
\begin{singlespace}
\noindent \label{prop:Nikolskii} Let $G\in L^{q}(\Omega_{T},\mathbb{R}^{k})$,
with $1\leq q<\infty$. Assume that for $\theta\in(0,1)$, $M\in[0,\infty)$
and a cylinder $Q_{\rho}(z_{0})=B_{\rho}(x_{0})\times(t_{0}-\rho^{2},t_{0})\Subset\Omega_{T}$,
the inequality 
\[
\left|h\right|^{-\,q\,\theta}\int_{Q_{\rho}(z_{0})}\left|G(x+he_{i},t)-G(x,t)\right|^{q}dx\,dt\,\leq\,M^{q}
\]
holds whenever $0<\left|h\right|\leq\min\left\{ \mathrm{dist}(B_{\rho}(x_{0}),\partial\Omega),A\right\} $
and $i\in\{1,\ldots,n\}$, where $A>0$ and $\{e_{i}\}_{1\,\leq\,i\,\leq\,n}$
is the standard basis of $\mathbb{R}^{n}$. Then 
\[
G\in L_{loc}^{q}\left(t_{0}-\rho^{2},t_{0};W_{loc}^{s,q}(B_{\rho}(x_{0}),\mathbb{R}^{k})\right)\,\,\,\,\,\,\,\,\mathit{for\,\,all}\,\,s\in(0,\theta).
\]
\end{singlespace}
\end{prop}

\section{The regularization \label{sec:a priori} }

\noindent $\hspace*{1em}$In this section, we shall assume that $\lambda$,
$\alpha$, $f$ and $u$ are as in the statement of Theorem \ref{thm:theo1}
if $p>2$, and as in the statement of Theorem \ref{thm:nuovo} if
$p=2$.

\noindent For $\varepsilon\in[0,1]$ and a couple of standard, non-negative,
radially symmetric mollifiers $\varphi_{1}\in C_{0}^{\infty}(B_{1}(0))$
and $\varphi_{2}\in C_{0}^{\infty}((-1,1))$, we define
\begin{equation}
f_{\varepsilon}(x,t):=\int_{-1}^{1}\int_{B_{1}(0)}f(x-\varepsilon y,t-\varepsilon s)\,\varphi_{1}(y)\,\varphi_{2}(s)\,dy\,ds,\label{eq:molli-1}
\end{equation}
where $f$ is meant to be extended by zero outside $\Omega_{T}$.
Let us observe that $f_{0}=f$ and $f_{\varepsilon}\in C^{\infty}(\Omega_{T})$
for every $\varepsilon\in(0,1]$. In addition, we define the function
$A_{\varepsilon}:\mathbb{R}^{n}\rightarrow\mathbb{R}$ by 
\[
A_{\varepsilon}(\xi)\,:=\,\frac{1}{p}\,(\vert\xi\vert-\lambda)_{+}^{p}\,+\,\frac{\varepsilon}{p}\,(1+\vert\xi\vert^{2})^{\frac{p}{2}}.
\]
Therefore, we obtain
\[
D_{\xi}A_{\varepsilon}(\xi)=\,H_{p-1}(\xi)+\varepsilon\,(1+\vert\xi\vert^{2})^{\frac{p-2}{2}}\xi\,,\,\,\,\,\,\,\,\,\xi\in\mathbb{R}^{n}.
\]
The next lemma provides the $p$-ellipticity property of the operator
$D_{\xi}A_{\varepsilon}(\xi)$ and its proof follows from \cite[Lemma 2.7]{BDGP}.
\selectlanguage{english}%
\begin{lem}
\noindent \label{lem:qform} Let $\varepsilon\in[0,1]$, $p\ge2$
and $z\in\mathbb{R}^{n}\setminus\{0\}$. Then, for every $\zeta\in\mathbb{R}^{n}$
we have
\[
\left[\varepsilon\,(1+|z|^{2})^{\frac{p-2}{2}}+\frac{(\vert z\vert-\lambda)_{+}^{p-1}}{\vert z\vert}\right]|\zeta|^{2}\le\langle D^{2}A_{\varepsilon}(z)\,\zeta,\zeta\rangle\le(p-1)\left[\varepsilon\,(1+|z|^{2})^{\frac{p-2}{2}}+(\vert z\vert-\lambda)_{+}^{p-2}\right]|\zeta|^{2}.
\]
\end{lem}

\selectlanguage{british}%
\noindent $\hspace*{1em}$Now we consider a parabolic cylinder $Q_{R}(z_{0}):=B_{R}(x_{0})\times(t_{0}-R^{2},t_{0})\Subset\Omega_{T}$.
For our purposes, in the following we will need the definition below.
\begin{defn}
\noindent \label{def:L2-traces} Let $\varepsilon\in(0,1]$. In this
framework, we identify a function 
\[
u_{\varepsilon}\in C^{0}\left([t_{0}-R^{2},t_{0}];L^{2}(B_{R}(x_{0}))\right)\cap L^{p}\left(t_{0}-R^{2},t_{0};W^{1,p}(B_{R}(x_{0}))\right)
\]
as a \textit{weak solution of the Cauchy-Dirichlet problem}
\begin{equation}
\begin{cases}
\begin{array}{cc}
\partial_{t}u_{\varepsilon}-\mathrm{div}\left[DA_{\varepsilon}(Du_{\varepsilon})\right]=f_{\varepsilon} & \mathrm{in}\,\,\,Q_{R}(z_{0}),\\
u_{\varepsilon}=u & \,\,\,\,\,\,\,\,\,\,\,\,\,\,\mathrm{on}\,\,\,\partial_{\mathrm{par}}Q_{R}(z_{0}),\,\,\,
\end{array}\end{cases}\label{eq:CAUCHYDIR}
\end{equation}
if and only if, for any $\varphi\in C_{0}^{\infty}(Q_{R}(z_{0}))$,
the following integral identity holds
\begin{equation}
\int_{Q_{R}(z_{0})}\left(u_{\varepsilon}\cdot\partial_{t}\varphi-\langle DA_{\varepsilon}(Du_{\varepsilon}),D\varphi\rangle\right)\,dz\,=\,-\int_{Q_{R}(z_{0})}f_{\varepsilon}\,\varphi\,dz\label{eq:weaksol2}
\end{equation}
and, moreover,
\[
u_{\varepsilon}\,\in\,u+L^{p}\left(t_{0}-R^{2},t_{0};W_{0}^{1,p}(B_{R}(x_{0}))\right)
\]
and $u_{\varepsilon}(\cdot,t_{0}-R^{2})=u(\cdot,t_{0}-R^{2})$ in
the $L^{2}$-sense, that is, 
\begin{equation}
\underset{t\,\downarrow\,t_{0}-R^{2}}{\lim}\,\big\Vert u_{\varepsilon}(\cdot,t)-u(\cdot,t_{0}-R^{2})\big\Vert_{L^{2}(B_{R}(x_{0}))}\,=\,0.\label{eq:L2sense}
\end{equation}
Therefore, the initial condition $u_{\varepsilon}=u$ on $B_{R}(x_{0})\times\{t_{0}-R^{2}\}$
has to be understood in the usual $L^{2}$-sense (\ref{eq:L2sense}),
while the condition $u_{\varepsilon}=u$ on the lateral boundary $\partial B_{R}(x_{0})\times(t_{0}-R^{2},t_{0})$
has to be meant in the sense of traces, i.e.
\[
(u_{\varepsilon}-u)\,(\cdot,t)\in W_{0}^{1,p}(B_{R}(x_{0}))
\]
for almost every $t\in(t_{0}-R^{2},t_{0})$.
\end{defn}

\noindent \begin{brem} The advantage of considering the Cauchy-Dirichlet
problem (\ref{eq:CAUCHYDIR}) stems from the fact that the existence
of a unique solution
\[
u_{\varepsilon}\,\in\,C^{0}\left([t_{0}-R^{2},t_{0}];L^{2}(B_{R}(x_{0}))\right)\cap L^{p}\left(t_{0}-R^{2},t_{0};W^{1,p}(B_{R}(x_{0}))\right)
\]
satisfying the requirements of Definition \ref{def:L2-traces} can
be ensured by the classic existence theory for parabolic equations
(see, for example, \cite[Chapter 2, Theorem 1.2 and Remark 1.2]{Lions}).
Moreover, the operator $D_{\xi}A_{\varepsilon}(\xi)$ fulfills $p$-growth
and $p$-ellipticity conditions with constants depending on $\varepsilon$
(cf. Lemma \ref{lem:qform} above). Therefore, by the results in \cite{Duzaar},
we have
\[
V_{1}(Du_{\varepsilon}):=(1+\vert Du_{\varepsilon}\vert^{2})^{\frac{p-2}{4}}Du_{\varepsilon}\,\in\,L_{loc}^{2}\left(t_{0}-R^{2},t_{0};W_{loc}^{1,2}(B_{R}(x_{0}),\mathbb{R}^{n})\right)
\]
and 
\[
Du_{\varepsilon}\,\in\,L_{loc}^{p+\frac{4}{n}}(Q_{R}(z_{0}),\mathbb{R}^{n}),
\]
and, by the definition of $V_{1}(Du_{\varepsilon})$, this yields
\[
DV_{1}(Du_{\varepsilon})\approx(1+\vert Du_{\varepsilon}\vert^{2})^{\frac{p-2}{4}}D^{2}u_{\varepsilon}\,\in\,L_{loc}^{2}(Q_{R}(z_{0}),\mathbb{R}^{n\times n})\,\,\Longrightarrow\,\,\vert D^{2}u_{\varepsilon}\vert\,\in\,L_{loc}^{2}(Q_{R}(z_{0})).
\]
\end{brem}

\subsection{Uniform a priori estimates}

\noindent $\hspace*{1em}$The first step in the proofs of Theorems
\ref{thm:theo1} and \ref{thm:nuovo} is the following estimate for
the weak solutions of the regularized problems (\ref{eq:CAUCHYDIR}).
Similar estimates are scattered in the literature, but not in the
exact form needed here. 
\begin{prop}[\textbf{Uniform energy estimate}]
\noindent  \label{prop:uniform1} With the notation and under the
assumptions above, there exist two positive constants $\varepsilon_{0}\leq1$
and $C\equiv C(n,p,R)$ such that 
\begin{equation}
\int_{Q_{R}(z_{0})}\vert Du_{\varepsilon}\vert^{p}\,dz\,\,+\sup_{t\in(t_{0}-R^{2},t_{0})}\big\Vert u_{\varepsilon}(\cdot,t)-u(\cdot,t)\big\Vert_{L^{2}(B_{R}(x_{0}))}^{2}\leq\,C\left(\Vert Du\Vert_{L^{p}(Q_{R}(z_{0}))}^{p}\,+\lambda^{p}+1\right)\label{eq:uni1}
\end{equation}
for all $\varepsilon\in(0,\varepsilon_{0}]$.
\end{prop}

\noindent \begin{proof}[\bfseries{Proof}] In order to obtain an (uniform
in $\varepsilon$) energy estimate for $\vert Du_{\varepsilon}\vert$,
we proceed by testing equations (\ref{eq:1}) and $(\ref{eq:CAUCHYDIR})_{1}$
with the map $\varphi=\psi(t)(u_{\varepsilon}-u)$, where $\psi\in W^{1,\infty}(\mathbb{R})$
is chosen such that
\[
\psi(t)=\begin{cases}
\begin{array}{c}
1\\
-\,\frac{1}{\omega}(t-t_{2}-\omega)\\
0
\end{array} & \begin{array}{c}
\mathrm{if}\,\,t\leq t_{2},\quad\quad\quad\;\,\,\,\,\\
\mathrm{if}\,\,t_{2}<t<t_{2}+\omega,\\
\mathrm{if}\,\,t\geq t_{2}+\omega,\quad\quad
\end{array}\end{cases}
\]
with
\[
t_{0}-R^{2}<t_{2}<t_{2}+\omega<t_{0},
\]
and then letting $\omega\rightarrow0$. We observe that, at this stage,
it is important that $u_{\varepsilon}$ and $u$ agree on the parabolic
boundary $\partial_{\mathrm{par}}Q_{R}(z_{0})$. We also note that
the following computations are somewhat formal concerning the use
of the time derivative, but they can easily be made rigorous, for
example by the use of Steklov averages. We skip this, since it is
a standard procedure. With the previous choice of $\varphi$, for
every $t_{2}\in(t_{0}-R^{2},t_{0})$ we find\begin{equation}\label{eq:est1}
\begin{split}
&\frac{1}{2}\int_{B_{R}(x_{0})}\left|u_{\varepsilon}(x,t_{2})-u(x,t_{2})\right|^{2}dx\,+\int_{Q_{R,t_{2}}}\langle H_{p-1}(Du_{\varepsilon})-H_{p-1}(Du),Du_{\varepsilon}-Du\rangle\,dz\\
& +\,\varepsilon\int_{Q_{R,t_{2}}}\langle(1+\vert Du_{\varepsilon}\vert^{2})^{\frac{p-2}{2}}Du_{\varepsilon},Du_{\varepsilon}-Du\rangle\,dz\,=\int_{Q_{R,t_{2}}}(f-f_{\varepsilon})(u_{\varepsilon}-u)\,dz,
\end{split}
\end{equation}where we have used the abbreviation
\[
Q_{R,t_{2}}=B_{R}(x_{0})\times(t_{0}-R^{2},t_{2}).
\]
In what follows, we shall denote by $c_{k}$ some positive constants
depending only on $n$, $p$ and $R$. Using Lemma \ref{lem:Brasco},
the Cauchy-Schwarz inequality as well as Young's inequality with $\beta>0$,
from \eqref{eq:est1} we infer\begin{equation}\label{eq:est2}
\begin{split}
&\sup_{t\in(t_{0}-R^{2},t_{0})}\big\Vert u_{\varepsilon}(\cdot,t)-u(\cdot,t)\big\Vert_{L^{2}(B_{R})}^{2}\,+\,\frac{8}{p^{2}}\int_{Q_{R}(z_{0})}\left|H_{\frac{p}{2}}(Du_{\varepsilon})-H_{\frac{p}{2}}(Du)\right|^{2}dz\,+\,2\,\varepsilon\int_{Q_{R}(z_{0})}\vert Du_{\varepsilon}\vert^{p}\,dz\\
&\,\,\,\,\,\,\,\leq\,2\,\varepsilon\int_{Q_{R}(z_{0})}(1+\vert Du_{\varepsilon}\vert^{2})^{\frac{p-1}{2}}\vert Du\vert\,dz\,+\,2\int_{Q_{R}(z_{0})}\left|f-f_{\varepsilon}\right|\left|u_{\varepsilon}-u\right|\,dz\\
&\,\,\,\,\,\,\,\leq\,\frac{2\,\varepsilon\beta^{p'}}{p'}\int_{Q_{R}}(1+\vert Du_{\varepsilon}\vert^{2})^{\frac{p}{2}}\,dz\,+\,\frac{2\,\varepsilon}{p\beta^{p}}\int_{Q_{R}}\vert Du\vert^{p}\,dz\,+\,2\int_{Q_{R}}\left|f-f_{\varepsilon}\right|\left|u_{\varepsilon}-u\right|\,dz\\
&\,\,\,\,\,\,\,\leq\,\frac{2^{\frac{p}{2}}\,\varepsilon\beta^{p'}}{p'}\int_{Q_{R}}\vert Du_{\varepsilon}\vert^{p}\,dz\,+\,\frac{2^{\frac{p}{2}}\,\varepsilon\beta^{p'}}{p'}\vert Q_{R}\vert+\,\frac{2\,\varepsilon}{p\beta^{p}}\int_{Q_{R}}\vert Du\vert^{p}\,dz\,+\,2\int_{Q_{R}}\left|f-f_{\varepsilon}\right|\left|u_{\varepsilon}-u\right|\,dz.
\end{split}
\end{equation}Choosing $\beta=\left(\frac{p'}{2^{p/2}}\right)^{\frac{1}{p'}}$ and
reabsorbing the first integral in the right-hand side of \eqref{eq:est2}
by the left-hand side, we arrive at\begin{equation}\label{eq:est3}
\begin{split}
&\sup_{t\in(t_{0}-R^{2},t_{0})}\big\Vert u_{\varepsilon}(\cdot,t)-u(\cdot,t)\big\Vert_{L^{2}(B_{R})}^{2}\,+\int_{Q_{R}}\left|H_{\frac{p}{2}}(Du_{\varepsilon})-H_{\frac{p}{2}}(Du)\right|^{2}dz\,+\,\varepsilon\int_{Q_{R}}\vert Du_{\varepsilon}\vert^{p}\,dz\\
&\,\,\,\,\,\,\,\leq\,c_{1}\,\varepsilon+\,c_{1}\,\varepsilon\int_{Q_{R}}\vert Du\vert^{p}\,dz\,+\,c_{1}\int_{Q_{R}}\left|f-f_{\varepsilon}\right|\left|u_{\varepsilon}-u\right|\,dz.
\end{split}
\end{equation}Now we apply Hölder's and Sobolev inequalities to estimate the last
integral as follows\begin{equation}\label{eq:est4}
\begin{split}
\int_{Q_{R}}\left|f-f_{\varepsilon}\right|\left|u_{\varepsilon}-u\right|\,dz&\,\leq\,\Vert f-f_{\varepsilon}\Vert_{L^{p'}(Q_{R})}\,\Vert u_{\varepsilon}-u\Vert_{L^{p}(Q_{R})}\\
&\leq \,c_{2}\,\Vert f-f_{\varepsilon}\Vert_{L^{p'}(Q_{R})}\,\Vert Du_{\varepsilon}-Du\Vert_{L^{p}(Q_{R})}.
\end{split}
\end{equation}Joining estimates \eqref{eq:est3} and \eqref{eq:est4}, recalling
that $0<\varepsilon\leq1$ and applying Young's inequality with $\gamma>0$,
we obtain\begin{align}\label{eq:est5}
&\sup_{t\in(t_{0}-R^{2},t_{0})}\big\Vert u_{\varepsilon}(\cdot,t)-u(\cdot,t)\big\Vert_{L^{2}(B_{R})}^{2}\,+\int_{Q_{R}}\left|H_{\frac{p}{2}}(Du_{\varepsilon})-H_{\frac{p}{2}}(Du)\right|^{2}dz\,+\,\varepsilon\int_{Q_{R}}\vert Du_{\varepsilon}\vert^{p}\,dz\nonumber\\
&\,\,\,\,\,\,\,\leq\,c_{3}\,+\,c_{3}\int_{Q_{R}}\vert Du\vert^{p}\,dz\,+\,c_{3}\,\gamma^{\frac{1}{1-p}}\Vert f-f_{\varepsilon}\Vert_{L^{p'}(Q_{R})}^{p'}+\,\gamma\,\Vert Du_{\varepsilon}-Du\Vert_{L^{p}(Q_{R})}^{p}\nonumber\\
&\,\,\,\,\,\,\,\leq\,c_{3}\,+(c_{3}+2^{p-1}\gamma)\int_{Q_{R}}\vert Du\vert^{p}\,dz\,+\,2^{p-1}\gamma\int_{Q_{R}}\vert Du_{\varepsilon}\vert^{p}\,dz\,+\,c_{3}\,\gamma^{\frac{1}{1-p}}\Vert f-f_{\varepsilon}\Vert_{L^{p'}(Q_{R})}^{p'}.
\end{align} Now, arguing as in \cite[Formula (4.20)]{AmPa}, we have
\begin{equation}
\int_{Q_{R}}\vert Du_{\varepsilon}\vert^{p}\,dz\,\leq\,2^{p}\int_{Q_{R}}\left|H_{\frac{p}{2}}(Du_{\varepsilon})-H_{\frac{p}{2}}(Du)\right|^{2}dz\,+\,2^{p+1}\int_{Q_{R}}(\vert Du\vert^{p}+\lambda^{p})\,dz.\label{eq:oldtrick}
\end{equation}
Multiplying all sides of \eqref{eq:est5} by $2^{p}$ and then adding
the resulting expression and (\ref{eq:oldtrick}) side by side, we
get\begin{align*}
&(1+2^{p}\varepsilon)\int_{Q_{R}}\vert Du_{\varepsilon}\vert^{p}\,dz\,+\,2^{p}\sup_{t\in(t_{0}-R^{2},t_{0})}\big\Vert u_{\varepsilon}(\cdot,t)-u(\cdot,t)\big\Vert_{L^{2}(B_{R})}^{2}\\
&\,\,\,\,\,\,\,\leq\,c_{4}+\,(c_{4}+2^{2p-1}\gamma)\int_{Q_{R}}(\vert Du\vert^{p}+\lambda^{p})\,dz\,+\,c_{4}\,\gamma^{\frac{1}{1-p}}\Vert f-f_{\varepsilon}\Vert_{L^{p'}(Q_{R})}^{p'}+\,2^{2p-1}\gamma\int_{Q_{R}}\vert Du_{\varepsilon}\vert^{p}\,dz.
\end{align*}Choosing $\gamma=\frac{1}{2^{2p}}$ and reabsorbing the last term
in the above estimate by the left-hand side, we find

\noindent \begin{equation}\label{eq:est6}
\begin{split}
&\int_{Q_{R}}\vert Du_{\varepsilon}\vert^{p}\,dz\,+\sup_{t\in(t_{0}-R^{2},t_{0})}\big\Vert u_{\varepsilon}(\cdot,t)-u(\cdot,t)\big\Vert_{L^{2}(B_{R})}^{2}\\
&\,\,\,\,\,\,\,\leq\,c_{5}\int_{Q_{R}}(\vert Du\vert^{p}+\lambda^{p})\,dz\,+\,c_{5}\,\Vert f-f_{\varepsilon}\Vert_{L^{p'}(Q_{R})}^{p'}+c_{5}\,.
\end{split}
\end{equation}Since, by virtue of (\ref{eq:molli-1}), 
\begin{equation}
f_{\varepsilon}\rightarrow f\,\,\,\,\,\,\,\,\,\mathrm{strongly\,\,in\,\,}L^{p'}(Q_{R})\,\,\,\,\,\mathrm{as}\,\,\varepsilon\rightarrow0^{+},\label{eq:strong}
\end{equation}
there exists a positive number $\varepsilon_{0}\leq1$ such that
\begin{equation}
\Vert f-f_{\varepsilon}\Vert_{L^{p'}(Q_{R})}\,\leq\,1\,\,\,\,\,\,\,\,\,\,\,\,\,\,\mathrm{for\,\,all}\,\,\varepsilon\in(0,\varepsilon_{0}].\label{eq:fconv}
\end{equation}
Combining this with \eqref{eq:est6}, we obtain the desired conclusion.\end{proof}

\noindent $\hspace*{1em}$Thanks to the previous proposition, we can
now obtain the following result, which will be decisive to carry out
the proofs of Theorems \ref{thm:theo1} and \ref{thm:nuovo}. 
\begin{prop}[\textbf{Comparison estimate}]
\noindent  \label{prop:comparison} With the notation and under the
assumptions of Proposition \ref{prop:uniform1}, there exists a positive
constant $C$ depending only on $n$, $p$ and $R$ such that the
estimate\begin{equation}\label{eq:comp}
\begin{split}
&\sup_{t\in(t_{0}-R^{2},t_{0})}\big\Vert u_{\varepsilon}(\cdot,t)-u(\cdot,t)\big\Vert_{L^{2}(B_{R}(x_{0}))}^{2}\,+\int_{Q_{R}(z_{0})}\left|\mathcal{V}_{\alpha,\lambda}(Du_{\varepsilon})-\mathcal{V}_{\alpha,\lambda}(Du)\right|^{2}dz\\
&\,\,\,\,\,\,\,\leq\,\,C\,\varepsilon\,\left(\Vert Du\Vert_{L^{p}(Q_{R}(z_{0}))}^{p}+1\right)+\,C\,\Vert f-f_{\varepsilon}\Vert_{L^{p'}(Q_{R}(z_{0}))}\,\left(\Vert Du\Vert_{L^{p}(Q_{R}(z_{0}))}+\lambda+1\right)
\end{split}
\end{equation} holds for every $\varepsilon\in(0,\varepsilon_{0}]$, where $\varepsilon_{0}$
is the constant from Proposition \ref{prop:uniform1}. In particular,
this estimate implies that
\[
u_{\varepsilon}\rightarrow u\,\,\,\,\,\,\,\,\,\mathit{strongly}\,\,\,\mathit{in}\,\,\,L^{2}(Q_{R}(z_{0}))\,\,\,\,\,\mathit{as}\,\,\,\varepsilon\rightarrow0
\]
and
\[
\mathcal{V}_{\alpha,\lambda}(Du_{\varepsilon})\rightarrow\mathcal{V}_{\alpha,\lambda}(Du)\,\,\,\,\,\,\,\,\,\mathit{strongly}\,\,\,\mathit{in}\,\,\,L^{2}(Q_{R}(z_{0}),\mathbb{R}^{n})\,\,\,\,\,\mathit{as}\,\,\,\varepsilon\rightarrow0.
\]
\end{prop}

\noindent \begin{proof}[\bfseries{Proof}] Arguing as in the proof
of Proposition \ref{prop:uniform1}, but using Lemma \ref{lem:vital}
instead of Lemma \ref{lem:Brasco}, we arrive at\begin{align*}
&\sup_{t\in(t_{0}-R^{2},t_{0})}\big\Vert u_{\varepsilon}(\cdot,t)-u(\cdot,t)\big\Vert_{L^{2}(B_{R})}^{2}\,+\int_{Q_{R}}\left|\mathcal{V}_{\alpha,\lambda}(Du_{\varepsilon})-\mathcal{V}_{\alpha,\lambda}(Du)\right|^{2}dz\\
&\,\,\,\,\,\,\,\leq\,c\,\varepsilon\,\left(\Vert Du\Vert_{L^{p}(Q_{R})}^{p}+1\right)\,+\,c\,\Vert f-f_{\varepsilon}\Vert_{L^{p'}(Q_{R})}\,\Vert Du_{\varepsilon}-Du\Vert_{L^{p}(Q_{R})}\,,
\end{align*}where $c\equiv c(n,p,R)>0$. Now, let us consider the same $\varepsilon_{0}\in(0,1]$
as in the proof of Proposition \ref{prop:uniform1}. Then, applying
Minkowski's inequality, for every $\varepsilon\in(0,\varepsilon_{0}]$
we get\begin{align*}
&\sup_{t\in(t_{0}-R^{2},t_{0})}\big\Vert u_{\varepsilon}(\cdot,t)-u(\cdot,t)\big\Vert_{L^{2}(B_{R})}^{2}\,+\int_{Q_{R}}\left|\mathcal{V}_{\alpha,\lambda}(Du_{\varepsilon})-\mathcal{V}_{\alpha,\lambda}(Du)\right|^{2}dz\\
&\,\,\,\,\,\,\,\leq\,c\,\varepsilon\,\left(\Vert Du\Vert_{L^{p}(Q_{R})}^{p}+1\right)+\,c\,\Vert f-f_{\varepsilon}\Vert_{L^{p'}(Q_{R})}\,\left(\Vert Du_{\varepsilon}\Vert_{L^{p}(Q_{R})}+\Vert Du\Vert_{L^{p}(Q_{R})}\right)\\
&\,\,\,\,\,\,\,\leq\,C\,\varepsilon\,\left(\Vert Du\Vert_{L^{p}(Q_{R})}^{p}+1\right)+\,C\,\Vert f-f_{\varepsilon}\Vert_{L^{p'}(Q_{R})}\,\left(\Vert Du\Vert_{L^{p}(Q_{R})}+\lambda+1\right),
\end{align*}where, in the last line, we have used inequality \eqref{eq:uni1}.
Thus we obtain the comparison estimate \eqref{eq:comp}. Now, just
observe that 
\[
\int_{Q_{R}(z_{0})}\vert u_{\varepsilon}-u\vert^{2}\,dz\,\,\leq\,\,R^{2}\,\sup_{t\in(t_{0}-R^{2},t_{0})}\big\Vert u_{\varepsilon}(\cdot,t)-u(\cdot,t)\big\Vert_{L^{2}(B_{R}(x_{0}))}^{2}.
\]
Combining this with \eqref{eq:comp} and (\ref{eq:strong}), the statement
is proved.\end{proof}

\noindent $\hspace*{1em}$We are now in a position to establish uniform
estimates for the derivatives of the function $\mathcal{V}_{\alpha,\lambda}(Du_{\varepsilon})$.
Let us start with the following result. 
\begin{prop}[\textbf{Uniform Sobolev estimate}]
\noindent  \label{prop:UniformSobolev} Under the assumptions of
Theorem \ref{thm:theo1} and with the notation above, there exists
a positive number $\varepsilon_{1}\leq1$ such that, for every $\varepsilon\in(0,\varepsilon_{1}]$,
for every parabolic cylinder $Q_{\rho}(z_{0})\Subset Q_{R}(z_{0})$
and every $j\in\{1,\ldots,n\}$, we have\begin{align}\label{eq:uniSobolev}
\int_{Q_{\rho/2}(z_{0})}\vert D_{x_{j}}\mathcal{V}_{\alpha,\lambda}(Du_{\varepsilon})\vert^{2}\,dz\,&\leq \left(C\,+\,\frac{C}{\rho^{2}}\right)\left[\Vert Du\Vert_{L^{p}(Q_{R})}^{p}+\Vert Du\Vert_{L^{p}(Q_{R})}^{2}+\lambda^{p}+\lambda^{2}+1\right]\nonumber\\
&\,\,\,\,\,\,\,+\,C\,\Vert f\Vert_{L^{p'}\left(t_{0}-R^{2},t_{0};B_{p',1}^{\frac{p-2}{p}}(B_{R}(x_{0}))\right)}^{p'}
\end{align}for a positive constant $C$ depending only on $n$, $p$ and $R$
in the case $\lambda=0$, and additionally on $\alpha$ if $\lambda>0$. 
\end{prop}

\noindent \begin{proof}[\bfseries{Proof}] Let us first assume that
$\lambda>0$. To shorten our notation, we introduce the function $P:Q_{R}(z_{0})\rightarrow\mathbb{R}_{0}^{+}$
defined by
\begin{equation}
P:=(\vert Du_{\varepsilon}\vert-\lambda)_{+}\label{eq:Pfunction}
\end{equation}
and its ``mollified'' version
\begin{equation}
P_{\varrho}\,:=\,(\vert[Du_{\varepsilon}]_{\varrho}\vert-\lambda)_{+}\,,\,\,\,\,\,\,\,\,\,\,\varrho>0,\label{eq:Pmolli}
\end{equation}
using the notation (\ref{eq:mollified}) on the right-hand side of
(\ref{eq:Pmolli}), while, on the left-hand side, an intentional abuse
of the same notation has been made. We now test the weak formulation
(\ref{eq:weaksol2}) with the map
\[
\varphi=D_{x_{j}}\left(D_{x_{j}}[u_{\varepsilon}]_{\varrho}\cdot\psi\cdot\Phi(P_{\varrho})\right),\,\,\,\,\,\,\,\,\,\,j\in\{1,\ldots,n\},
\]
where $\psi\in W_{0}^{1,\infty}(Q_{R}(z_{0}))$ is a non-negative
cut-off function that will be specified later and\foreignlanguage{english}{
$\Phi:\mathbb{R}_{0}^{+}\rightarrow\mathbb{R}_{0}^{+}$ is an increasing,
locally Lipschitz continuous function, such that $\Phi$ and $\Phi'$
are bounded on $\mathbb{R}_{0}^{+}$, $\Phi(0)=0$ and 
\begin{equation}
\Phi'(t)\,t\,\leq\,c_{\Phi}\,\Phi(t)\label{eq:derivative}
\end{equation}
for a suitable constant $c_{\Phi}>0$.} After summing with respect
to $j$ from $1$ to $n$, we obtain\begin{equation}\label{eq:new1}
\begin{split}
&\frac{1}{2}\int_{Q_{R}}\partial_{t}(\vert D[u_{\varepsilon}]_{\varrho}\vert^{2})\,\psi\,\Phi(P_{\varrho})\,dz\,+\int_{Q_{R}}\sum_{j=1}^{n}\left\langle D_{x_{j}}[DA_{\varepsilon}(Du_{\varepsilon})]_{\varrho},D\left(D_{x_{j}}[u_{\varepsilon}]_{\varrho}\cdot\psi\,\Phi(P_{\varrho})\right)\right\rangle dz\\
&\,\,\,\,\,\,\,= \,-\int_{Q_{R}}\sum_{j=1}^{n}\,[f_{\varepsilon}]_{\varrho}\,D_{x_{j}}\left(D_{x_{j}}[u_{\varepsilon}]_{\varrho}\cdot\psi\,\Phi(P_{\varrho})\right)dz.
\end{split}
\end{equation}Arguing as in \cite[page 521]{AmBa}, the integral involving the time
derivative can be expressed as follows\begin{align*}
\frac{1}{2}\int_{Q_{R}}\partial_{t}(\vert D[u_{\varepsilon}]_{\varrho}\vert^{2})\,\psi\,\Phi(P_{\varrho})\,dz\,&=\,\frac{1}{2}\int_{Q_{R}}\partial_{t}\left[\int_{0}^{\vert D[u_{\varepsilon}]_{\varrho}\vert^{2}}\Phi((\sqrt{w}-\lambda)_{+})\,dw\right]\psi\,dz\\
&=\,-\,\frac{1}{2}\int_{Q_{R}}\left[\int_{0}^{\vert D[u_{\varepsilon}]_{\varrho}\vert^{2}}\Phi((\sqrt{w}-\lambda)_{+})\,dw\right]\partial_{t}\psi\,dz.
\end{align*}Therefore, integrating by parts and then letting $\varrho\rightarrow0$
in \eqref{eq:new1}, we get\begin{align}\label{eq:new2}
&-\,\frac{1}{2}\int_{Q_{R}}\left[\int_{0}^{\vert Du_{\varepsilon}\vert^{2}}\Phi((\sqrt{w}-\lambda)_{+})\,dw\right]\partial_{t}\psi\,dz\nonumber\\
&+\int_{Q_{R}}\sum_{j=1}^{n}\,\langle D^{2}A_{\varepsilon}(Du_{\varepsilon})\,D(D_{x_{j}}u_{\varepsilon}),D(D_{x_{j}}u_{\varepsilon})\rangle \,\psi\,\Phi(P)\,dz\nonumber\\
&+\int_{Q_{R}}\sum_{j=1}^{n}\,\langle D^{2}A_{\varepsilon}(Du_{\varepsilon})\,D(D_{x_{j}}u_{\varepsilon}),D\psi\rangle (D_{x_{j}}u_{\varepsilon})\,\Phi(P)\,dz\nonumber\\
&+\int_{Q_{R}}\sum_{j=1}^{n}\,\langle D^{2}A_{\varepsilon}(Du_{\varepsilon})\,D(D_{x_{j}}u_{\varepsilon}),DP\rangle (D_{x_{j}}u_{\varepsilon})\,\psi\,\Phi'(P)\,dz\nonumber\\
&\,\,\,\,\,\,\,=\int_{Q_{R}}\sum_{j=1}^{n}\,(D_{x_{j}}f_{\varepsilon})\,(D_{x_{j}}u_{\varepsilon})\,\psi\,\Phi(P)\,dz.
\end{align}$\hspace*{1em}$In what follows, we will denote by $c_{k}$ and $C$
some positive constants that do not depend on $\varepsilon$. Now,
let us fix an arbitrary radius $\rho\in(0,R)$. For a fixed time $t_{1}\in(t_{0}-\rho^{2},t_{0})$
and $\sigma\in(0,t_{0}-t_{1})$, we choose
\[
\psi(x,t)=\widetilde{\chi}(t)\,\chi(t)\,\eta^{2}(x)
\]
with $\chi\in W^{1,\infty}((t_{0}-R^{2},t_{0}),[0,1])$, $\chi\equiv0$
on $(t_{0}-R^{2},t_{0}-\rho^{2})$ and $\partial_{t}\chi\geq0$, $\eta\in C_{0}^{\infty}(B_{\rho}(x_{0}),[0,1])$,
and with the Lipschitz continuous function $\widetilde{\chi}:(t_{0}-R^{2},t_{0})\rightarrow\mathbb{R}$
defined by
\[
\widetilde{\chi}(t)=\begin{cases}
\begin{array}{c}
1\\
\mathrm{affine}\\
0
\end{array} & \begin{array}{c}
\mathrm{if}\,\,t\leq t_{1},\quad\quad\quad\;\,\,\,\,\\
\mathrm{if}\,\,t_{1}<t<t_{1}+\sigma,\\
\mathrm{if}\,\,t\geq t_{1}+\sigma.\quad\quad
\end{array}\end{cases}
\]

\noindent With such a choice of $\psi$, equation \eqref{eq:new2}
turns into\begin{align*}
&-\,\frac{1}{2}\int_{Q_{R}}\left[\int_{0}^{\vert Du_{\varepsilon}\vert^{2}}\Phi((\sqrt{w}-\lambda)_{+})\,dw\right]\eta^{2}(x)\,\chi(t)\,\partial_{t}\widetilde{\chi}(t)\,dz\\
&-\,\frac{1}{2}\int_{Q_{R}}\left[\int_{0}^{\vert Du_{\varepsilon}\vert^{2}}\Phi((\sqrt{w}-\lambda)_{+})\,dw\right]\eta^{2}(x)\,\widetilde{\chi}(t)\,\partial_{t}\chi(t)\,dz\\
&+\int_{Q_{R}}\sum_{j=1}^{n}\,\langle D^{2}A_{\varepsilon}(Du_{\varepsilon})\,D(D_{x_{j}}u_{\varepsilon}),D(D_{x_{j}}u_{\varepsilon})\rangle \,\eta^{2}(x)\,\widetilde{\chi}(t)\,\chi(t)\,\Phi(P)\,dz\\
&+\int_{Q_{R}}\sum_{j=1}^{n}\,\langle D^{2}A_{\varepsilon}(Du_{\varepsilon})\,D(D_{x_{j}}u_{\varepsilon}),DP\rangle (D_{x_{j}}u_{\varepsilon})\,\eta^{2}(x)\,\widetilde{\chi}(t)\,\chi(t)\,\Phi'(P)\,dz\\
&\,\,\,\,\,\,\,=\,-\,2\int_{Q_{R}}\sum_{j=1}^{n}\,\langle D^{2}A_{\varepsilon}(Du_{\varepsilon})\,D(D_{x_{j}}u_{\varepsilon}),D\eta\rangle (D_{x_{j}}u_{\varepsilon})\,\eta(x)\,\widetilde{\chi}(t)\,\chi(t)\,\Phi(P)\,dz\\
&\,\,\,\,\,\,\,\,\,\,\,\,\,\,+\int_{Q_{R}}\sum_{j=1}^{n}\,(D_{x_{j}}f_{\varepsilon})\,(D_{x_{j}}u_{\varepsilon})\,\eta^{2}(x)\,\widetilde{\chi}(t)\,\chi(t)\,\Phi(P)\,dz.
\end{align*}Setting
\[
Q^{t_{1}}:=B_{\rho}(x_{0})\times(t_{0}-\rho^{2},t_{1})
\]
and letting $\sigma\rightarrow0$ in the previous equality, for every
$t_{1}\in(t_{0}-\rho^{2},t_{0})$ we get\begin{align}\label{eq:new3}
&\frac{1}{2}\int_{B_{\rho}(x_{0})}\chi(t_{1})\,\eta^{2}(x)\left[\int_{0}^{\vert Du_{\varepsilon}(x,t_{1})\vert^{2}}\Phi((\sqrt{w}-\lambda)_{+})\,dw\right]dx\nonumber\\
&+\int_{Q^{t_{1}}}\sum_{j=1}^{n}\,\langle D^{2}A_{\varepsilon}(Du_{\varepsilon})\,D(D_{x_{j}}u_{\varepsilon}),D(D_{x_{j}}u_{\varepsilon})\rangle \,\eta^{2}(x)\,\chi(t)\,\Phi(P)\,dz\nonumber\\
&+\int_{Q^{t_{1}}}\sum_{j=1}^{n}\,\langle D^{2}A_{\varepsilon}(Du_{\varepsilon})\,D(D_{x_{j}}u_{\varepsilon}),DP\rangle (D_{x_{j}}u_{\varepsilon})\,\eta^{2}(x)\,\chi(t)\,\Phi'(P)\,dz\nonumber\\
&\,\,\,\,\,\,\,=\,-\,2\int_{Q^{t_{1}}}\sum_{j=1}^{n}\,\langle D^{2}A_{\varepsilon}(Du_{\varepsilon})\,D(D_{x_{j}}u_{\varepsilon}),D\eta\rangle (D_{x_{j}}u_{\varepsilon})\,\eta(x)\,\chi(t)\,\Phi(P)\,dz\nonumber\\
&\,\,\,\,\,\,\,\,\,\,\,\,\,\,\,+\int_{Q^{t_{1}}}\sum_{j=1}^{n}\,(D_{x_{j}}f_{\varepsilon})\,(D_{x_{j}}u_{\varepsilon})\,\eta^{2}(x)\,\chi(t)\,\Phi(P)\,dz\nonumber\\
&\,\,\,\,\,\,\,\,\,\,\,\,\,\,+\,\frac{1}{2}\int_{Q^{t_{1}}}(\partial_{t}\chi)\,\eta^{2}(x)\left[\int_{0}^{\vert Du_{\varepsilon}(x,t)\vert^{2}}\Phi((\sqrt{w}-\lambda)_{+})\,dw\right]dz,
\end{align}where we have used that $\partial_{t}\widetilde{\chi}$ converges
to a Dirac delta distribution as $\sigma\rightarrow0$, together with
the $L^{2}(B_{R}(x_{0}))$-valued continuity of $u_{\varepsilon}$.
We now observe that the second integral on the left-hand side of \eqref{eq:new3}
is non-negative, since $A_{\varepsilon}$ is convex and $\chi,\Phi\geq0$.
As for the first term on the right-hand side of \eqref{eq:new3},
for every $t_{1}\in(t_{0}-\rho^{2},t_{0})$ we have\begin{align}\label{eq:new4}
&-\,2\int_{Q^{t_{1}}}\langle D^{2}A_{\varepsilon}(Du_{\varepsilon})\,D(D_{x_{j}}u_{\varepsilon}),D\eta\rangle (D_{x_{j}}u_{\varepsilon})\,\eta(x)\,\chi(t)\,\Phi(P)\,dz\nonumber\\
&\,\,\,\leq2\int_{Q^{t_{1}}}\vert\langle D^{2}A_{\varepsilon}(Du_{\varepsilon})\,D(D_{x_{j}}u_{\varepsilon}),D\eta\rangle\vert\,\chi(t)\,\eta(x)\,\vert D_{x_{j}}u_{\varepsilon}\vert\,\Phi (P)\,dz\nonumber\\
&\,\,\,\leq\,2\int_{Q^{t_{1}}}\sqrt{\langle D^{2}A_{\varepsilon}(Du_{\varepsilon})\,D(D_{x_{j}}u_{\varepsilon}),D(D_{x_{j}}u_{\varepsilon})\rangle}\,\sqrt{\langle D^{2}A_{\varepsilon}(Du_{\varepsilon})\,D\eta,D\eta\rangle}\,\chi(t)\,\eta(x)\,\vert D_{x_{j}}u_{\varepsilon}\vert\,\Phi(P)\,dz\nonumber\\
&\,\,\,\leq\,\frac{1}{2}\int_{Q^{t_{1}}}\langle D^{2}A_{\varepsilon}(Du_{\varepsilon})\,D(D_{x_{j}}u_{\varepsilon}),D(D_{x_{j}}u_{\varepsilon})\rangle\,\chi(t)\,\eta^{2}(x)\,\Phi (P)\,dz\nonumber\\
&\,\,\,\,\,\,\,\,\,\,+\,\,2\int_{Q^{t_{1}}}\langle D^{2}A_{\varepsilon}(Du_{\varepsilon})\,D\eta,D\eta\rangle\,\chi(t)\,\vert D_{x_{j}}u_{\varepsilon}\vert^{2}\,\Phi (P)\,dz,
\end{align}where we have used the Cauchy-Schwarz inequality together with Young's
inequality. Joining \eqref{eq:new3} and \eqref{eq:new4}, we obtain
\begin{equation}
I_{1}+I_{2}+I_{3}\leq\,I_{4}+I_{5}+I_{6}\,,\label{eq:integrals}
\end{equation}
where\begin{align*}
&I_{1}:=\int_{B_{\rho}(x_{0})}\chi(t_{1})\,\eta^{2}(x)\left[\int_{0}^{\vert Du_{\varepsilon}(x,t_{1})\vert^{2}}\Phi((\sqrt{w}-\lambda)_{+})\,dw\right]dx,\\
&I_{2}:=\,2\int_{Q^{t_{1}}}\sum_{j=1}^{n}\,\langle D^{2}A_{\varepsilon}(Du_{\varepsilon})\,D(D_{x_{j}}u_{\varepsilon}),DP\rangle(D_{x_{j}}u_{\varepsilon})\,\eta^{2}(x)\,\chi(t)\,\Phi'(P)\,dz,\\
&I_{3}:=\int_{Q^{t_{1}}}\sum_{j=1}^{n}\,\langle D^{2}A_{\varepsilon}(Du_{\varepsilon})\,D(D_{x_{j}}u_{\varepsilon}),D(D_{x_{j}}u_{\varepsilon})\rangle\,\chi(t)\,\eta^{2}(x)\,\Phi(P)\,dz,\\
&I_{4}:=\,4\int_{Q^{t_{1}}}\sum_{j=1}^{n}\,\langle D^{2}A_{\varepsilon}(Du_{\varepsilon})\,D\eta,D\eta\rangle\,\chi(t)\,\vert D_{x_{j}}u_{\varepsilon}\vert^{2}\,\Phi(P)\,dz,\\
&I_{5}:=\int_{Q^{t_{1}}}(\partial_{t}\chi)\,\eta^{2}(x)\left[\int_{0}^{\vert Du_{\varepsilon}(x,t)\vert^{2}}\Phi((\sqrt{w}-\lambda)_{+})\,dw\right]dz,\\
&I_{6}:=\,2\int_{Q^{t_{1}}}\sum_{j=1}^{n}\,(D_{x_{j}}f_{\varepsilon})\,(D_{x_{j}}u_{\varepsilon})\,\eta^{2}(x)\,\chi(t)\,\Phi(P)\,dz.
\end{align*}Now we observe that $I_{1}\geq0$, thus we can drop it in the following.
Furthermore, arguing as in the proof of \cite[Proposition 5.1]{AmGrPa},
one can show that $I_{2}\geq0$. Therefore, inequality (\ref{eq:integrals})
boils down to
\begin{equation}
I_{3}\leq\,I_{4}+I_{5}+I_{6}\,.\label{eq:integrals2}
\end{equation}
$\hspace*{1em}$At this stage, we choose a cut-off function $\eta\in C_{0}^{\infty}(B_{\rho}(x_{0}))$
with $\eta\equiv1$ on $B_{\rho/2}(x_{0})$ such that
\begin{equation}
0\leq\eta\leq1\,\,\,\,\,\,\,\,\,\,\mathrm{and}\,\,\,\,\,\,\,\,\,\,\left|D\eta\right|\leq\,\frac{\tilde{c}}{\rho}\,.\label{eq:proeta}
\end{equation}
For the cut-off function in time, we choose the piecewise affine function
$\chi:(t_{0}-R^{2},t_{0})\rightarrow\left[0,1\right]$ with
\[
\chi\equiv0\,\,\,\,\,\,\,\,\mathrm{on}\,\,\,(t_{0}-R^{2},t_{0}-\rho^{2}),\,
\]
\[
\chi\equiv1\,\,\,\,\,\,\,\,\mathrm{on}\,\,\,\left(t_{0}-\left(\frac{\rho}{2}\right)^{2},t_{0}\right),
\]
\[
\partial_{t}\chi\equiv\frac{4}{3\rho^{2}}\,\,\,\,\,\,\,\,\mathrm{on}\,\,\,\left(t_{0}-\rho^{2},t_{0}-\left(\frac{\rho}{2}\right)^{2}\right).
\]
Moreover, \foreignlanguage{english}{as in \cite{AmGrPa} we choose
\begin{equation}
\Phi(t)\,:=\,\frac{t^{2\alpha}}{(t^{2}+\lambda^{2})^{\alpha}}\,\,\,\,\,\,\,\,\,\,\mathrm{for}\,\,t\geq0,\label{eq:Phi}
\end{equation}
and therefore 
\begin{equation}
\Phi'(t)\,=\,\frac{2\alpha\lambda^{2}\,t^{2\alpha-1}}{(t^{2}+\lambda^{2})^{\alpha+1}}\,.\label{eq:derPhi}
\end{equation}
This function satisfies (\ref{eq:derivative}) with $c_{\Phi}=2\alpha$.}
With the above choices, we now estimate $I_{3}$, $I_{4}$ and $I_{5}$
separately. Let us first consider $I_{3}$. Recalling the definition
of $P$ in (\ref{eq:Pfunction}), by Lemma \ref{lem:qform} we have
\begin{equation}
I_{3}\,\geq\int_{Q^{t_{1}}}\frac{(\vert Du_{\varepsilon}\vert-\lambda)_{+}^{p-1+2\alpha}}{\vert Du_{\varepsilon}\vert\,[\lambda^{2}+(\vert Du_{\varepsilon}\vert-\lambda)_{+}^{2}]^{\alpha}}\,\vert D^{2}u_{\varepsilon}\vert^{2}\,\chi(t)\,\eta^{2}(x)\,dz.\label{eq:I3}
\end{equation}
Using Lemma \ref{lem:qform} again, the fact that $\varepsilon,\chi,\Phi\leq1$
and the properties (\ref{eq:proeta}) of $\eta$, we obtain
\begin{equation}
I_{4}\,\leq\,8\,(p-1)\int_{Q_{R}}(1+\vert Du_{\varepsilon}\vert^{2})^{\frac{p}{2}}\,\vert D\eta\vert^{2}\,dz\,\leq\,\frac{c_{1}(p)}{\rho^{2}}\int_{Q_{R}}(1+\vert Du_{\varepsilon}\vert^{2})^{\frac{p}{2}}\,dz.\label{eq:I4}
\end{equation}
Now we use again the fact that $\Phi\leq1$, the properties of $\chi$
and $\eta$, and Hölder's inequality, in order to get 
\begin{equation}
I_{5}\,\leq\int_{Q^{t_{1}}}(\partial_{t}\chi)\,\eta^{2}(x)\,\vert Du_{\varepsilon}\vert^{2}\,dz\,\leq\,\frac{4}{3\rho^{2}}\int_{Q_{R}}\vert Du_{\varepsilon}\vert^{2}\,dz\,\le\,\frac{c_{2}}{\rho^{2}}\,\Vert Du_{\varepsilon}\Vert_{L^{p}(Q_{R})}^{2}\,,\label{eq:I5}
\end{equation}
where $c_{2}\equiv c_{2}(n,p,R)>0$. Combining estimates (\ref{eq:integrals2}),
(\ref{eq:I3}), (\ref{eq:I4}) and (\ref{eq:I5}), we find\begin{align}\label{eq:key1}
&\int_{Q^{t_{1}}}\frac{(\vert Du_{\varepsilon}\vert-\lambda)_{+}^{p-1+2\alpha}}{\vert Du_{\varepsilon}\vert\,[\lambda^{2}+(\vert Du_{\varepsilon}\vert-\lambda)_{+}^{2}]^{\alpha}}\,\vert D^{2}u_{\varepsilon}\vert^{2}\,\chi(t)\,\eta^{2}(x)\,dz\nonumber\\
&\,\,\,\,\,\,\,\leq\,\frac{c_{2}}{\rho^{2}}\left[\int_{Q_{R}}(1+\vert Du_{\varepsilon}\vert^{2})^{\frac{p}{2}}\,dz\,+\,\Vert Du_{\varepsilon}\Vert_{L^{p}(Q_{R})}^{2}\right]\nonumber\\
&\,\,\,\,\,\,\,\,\,\,\,\,\,\,+\,2\,\sum_{j=1}^{n}\int_{Q^{t_{1}}}(D_{x_{j}}f_{\varepsilon})\,(D_{x_{j}}u_{\varepsilon})\,\eta^{2}(x)\,\chi(t)\,\Phi(P)\,dz,
\end{align}which holds for every $t_{1}\in(t_{0}-\rho^{2},t_{0})$.\\
$\hspace*{1em}$It remains to estimate the integral containing $D_{x_{j}}f_{\varepsilon}$.
With this aim in mind, we now argue as in the proof of \cite[Proposition 5.1]{AmGrPa}.
\foreignlanguage{english}{By Theorem \ref{thm:extension}, there exists
a linear and bounded extension operator 
\[
\mathrm{ext}_{B_{\rho}(x_{0})}:B_{p',1}^{-2/p}(B_{\rho}(x_{0}))\hookrightarrow B_{p',1}^{-2/p}(\mathbb{R}^{n})
\]
such that $\mathrm{re}_{B_{\rho}(x_{0})}\circ\mathrm{ext}_{B_{\rho}(x_{0})}=\mathrm{id}$,
where $\mathrm{re}_{B_{\rho}(x_{0})}$ is the restriction operator
defined in Section \ref{sec:Functional spaces} and the symbol $\mathrm{id}$
denotes the identity in $B_{p',1}^{-2/p}(B_{\rho}(x_{0}))$. Since,
for almost every $t\in(t_{0}-\rho^{2},t_{0})$, we have $D_{x_{j}}f_{\varepsilon}(\cdot,t)=\mathrm{ext}_{B_{\rho}(x_{0})}(D_{x_{j}}f_{\varepsilon}(\cdot,t))$
almost everywhere in $B_{\rho}(x_{0})$, we find}
\begin{equation}
\int_{Q^{t_{1}}}(D_{x_{j}}f_{\varepsilon})(D_{x_{j}}u_{\varepsilon})\,\eta^{2}(x)\,\chi(t)\,\Phi(P)\,dz\,=\int_{Q^{t_{1}}}\mathrm{ext}_{B_{\rho}(x_{0})}(D_{x_{j}}f_{\varepsilon}(x,t))\cdot(D_{x_{j}}u_{\varepsilon})\,\eta^{2}(x)\,\chi(t)\,\Phi(P)\,dz.\label{eq:pivot00}
\end{equation}
By definition of dual norm, we get\begin{equation}\label{eq:pivot}
\begin{split}
&\left|\int_{Q^{t_{1}}}\mathrm{ext}_{B_{\rho}(x_{0})}(D_{x_{j}}f_{\varepsilon}(x,t))\cdot(D_{x_{j}}u_{\varepsilon})\,\eta^{2}(x)\,\chi(t)\,\Phi(P)\,dz\right|\\
&\,\,\,\leq\int_{t_{0}-\rho^{2}}^{t_{1}}\Vert\mathrm{ext}_{B_{\rho}(x_{0})}(D_{x_{j}}f_{\varepsilon}(\cdot,t))\Vert_{(\mathring{B}_{p,\infty}^{2/p}(\mathbb{R}^{n}))'}\,\Vert(D_{x_{j}}u_{\varepsilon}(\cdot,t))\,\eta^{2}\,\chi(t)\,\Phi(P(\cdot,t))\Vert_{B_{p,\infty}^{2/p}(\mathbb{R}^{n})}\,dt\\
&\,\,\,=\int_{t_{0}-\rho^{2}}^{t_{1}}\Vert\mathrm{ext}_{B_{\rho}(x_{0})}(D_{x_{j}}f_{\varepsilon}(\cdot,t))\Vert_{(\mathring{B}_{p,\infty}^{2/p}(\mathbb{R}^{n}))'}\\
&\,\,\,\,\,\,\,\,\,\,\,\,\,\,\,\,\,\,\,\,\,\,\cdot\left(\Vert(D_{x_{j}}u_{\varepsilon}(\cdot,t))\,\eta^{2}\,\chi(t)\,\Phi(P(\cdot,t))\Vert_{L^{p}(\mathbb{R}^{n})}\,+\,[(D_{x_{j}}u_{\varepsilon}(\cdot,t))\,\eta^{2}\,\chi(t)\,\Phi(P(\cdot,t))]_{B_{p,\infty}^{2/p}(\mathbb{R}^{n})}\right)dt,
\end{split}
\end{equation}\foreignlanguage{english}{where, in the last line, we have used the
equivalence (\ref{eq:equivalence}).} By the properties of $\eta$
and the fact that $0\leq\chi,\Phi\leq1$, we have \foreignlanguage{english}{
\begin{equation}
\Vert(D_{x_{j}}u_{\varepsilon}(\cdot,t))\,\eta^{2}\,\chi(t)\,\Phi(P(\cdot,t))\Vert_{L^{p}(\mathbb{R}^{n})}\,\leq\,\Vert Du_{\varepsilon}(\cdot,t)\Vert_{L^{p}(B_{\rho}(x_{0}))}\label{eq:extra}
\end{equation}
for almost every $t\in(t_{0}-\rho^{2},t_{0})$. Furthermore, using
Theorem \ref{duality}, for almost every $t\in(t_{0}-\rho^{2},t_{0})$
we obtain 
\[
\Vert\mathrm{ext}_{B_{\rho}(x_{0})}(D_{x_{j}}f_{\varepsilon}(\cdot,t))\Vert_{(\mathring{B}_{p,\infty}^{2/p}(\mathbb{R}^{n}))'}\,\le\,c_{3}\,\Vert\mathrm{ext}_{B_{\rho}(x_{0})}(D_{x_{j}}f_{\varepsilon}(\cdot,t))\Vert_{B_{p',1}^{-2/p}(\mathbb{R}^{n})}\,,
\]
for some positive constant $c_{3}\equiv c_{3}(n,p)$. Combining the
above inequality with the boundedness of the operator $\mathrm{ext}_{B_{\rho}(x_{0})}$
yields 
\[
\Vert\mathrm{ext}_{B_{\rho}(x_{0})}(D_{x_{j}}f_{\varepsilon}(\cdot,t))\Vert_{(\mathring{B}_{p,\infty}^{2/p}(\mathbb{R}^{n}))'}\,\leq\,c_{3}\,\Vert D_{x_{j}}f_{\varepsilon}(\cdot,t)\Vert_{B_{p',1}^{-2/p}(B_{\rho}(x_{0}))}
\]
for almost every $t\in(t_{0}-\rho^{2},t_{0})$. Moreover, applying
Theorem \ref{negder}, we find that 
\[
\Vert D_{x_{j}}f_{\varepsilon}(\cdot,t)\Vert_{B_{p',1}^{-2/p}(B_{\rho}(x_{0}))}\,\le\,c_{3}\,\Vert f_{\varepsilon}(\cdot,t)\Vert_{B_{p',1}^{\frac{p-2}{p}}(B_{\rho}(x_{0}))}
\]
for almost every $t\in(t_{0}-\rho^{2},t_{0})$. Combining the preceding
inequalities, we infer 
\begin{equation}
\Vert\mathrm{ext}_{B_{\rho}(x_{0})}(D_{x_{j}}f_{\varepsilon}(\cdot,t))\Vert_{(\mathring{B}_{p,\infty}^{2/p}(\mathbb{R}^{n}))'}\,\le\,c_{3}\,\Vert f_{\varepsilon}(\cdot,t)\Vert_{B_{p',1}^{\frac{p-2}{p}}(B_{\rho}(x_{0}))},\label{eq:derf}
\end{equation}
which holds for almost every $t\in(t_{0}-\rho^{2},t_{0})$ and a positive
constant $c_{3}$ depending only on $n$ and $p$. Now, recalling
that\begin{align*}
&[(D_{x_{j}}u_{\varepsilon}(\cdot,t))\,\eta^{2}\,\chi(t)\,\Phi(P(\cdot,t))]_{B_{p,\infty}^{2/p}(\mathbb{R}^{n})}^{p}\\
&\,\,\,\,\,\,\,=\,\sup_{\vert h\vert\,>\,0}\,\frac{1}{\vert h\vert^{2}}\int_{\mathbb{R}^{n}}\chi^{p}(t)\left|\delta_{h}\left((D_{x_{j}}u_{\varepsilon}(x,t))\,\eta^{2}(x)\,\Phi(P(x,t))\right)\right|^{p}dx
\end{align*} and using Lemma \ref{lem:Lind}, we deduce
\begin{align*}
 & \frac{1}{\vert h\vert^{2}}\int_{\mathbb{R}^{n}}\chi^{p}(t)\left|\delta_{h}\left((D_{x_{j}}u_{\varepsilon}(x,t))\,\eta^{2}(x)\,\Phi(P(x,t))\right)\right|^{p}dx\\
 & \,\,\,\,\,\,\,\leq\,\frac{c_{4}(p)}{\vert h\vert^{2}}\,\chi(t)\int_{\mathbb{R}^{n}}\left|\delta_{h}\left(\vert D_{x_{j}}u_{\varepsilon}(x,t)\vert^{\frac{p-2}{2}}(D_{x_{j}}u_{\varepsilon}(x,t))\,\eta^{p}(x)\left[\Phi(P(x,t))\right]^{\frac{p}{2}}\right)\right|^{2}dx\\
 & \,\,\,\,\,\,\,\leq c_{5}(n,p)\,\chi(t)\int_{\mathbb{R}^{n}}\left|D\left(\vert D_{x_{j}}u_{\varepsilon}(x,t)\vert^{\frac{p-2}{2}}(D_{x_{j}}u_{\varepsilon}(x,t))\,\eta^{p}(x)\left[\Phi(P(x,t))\right]^{\frac{p}{2}}\right)\right|^{2}dx,
\end{align*}
where, in the last line, we have used the first statement in Lemma
\ref{lem:Giusti1}. Using the properties (\ref{eq:proeta}) of $\eta$
and the boundedness of $\Phi$, with simple manipulations we thus
obtain\begin{align}\label{eq:new5}
&[(D_{x_{j}}u_{\varepsilon}(\cdot,t))\,\eta^{2}\,\chi(t)\,\Phi(P(\cdot,t))]_{B_{p,\infty}^{2/p}(\mathbb{R}^{n})}^{p}\nonumber\\
&\,\,\,\,\,\,\,\leq\,c_{5}\,\chi(t)\int_{B_{\rho}(x_{0})}\left|D\left(\vert D_{x_{j}}u_{\varepsilon}(x,t)\vert^{\frac{p-2}{2}}(D_{x_{j}}u_{\varepsilon}(x,t))\left[\Phi(P(x,t))\right]^{\frac{p}{2}}\right)\right|^{2}\eta^{2}(x)\,dx\nonumber\\
&\,\,\,\,\,\,\,\,\,\,\,\,\,\,+\,\frac{c_{5}}{\rho^{2}}\,\chi(t)\int_{B_{\rho}(x_{0})}\vert Du_{\varepsilon}(x,t)\vert^{p}\,dx
\end{align} }for almost every $t\in(t_{0}-\rho^{2},t_{1})$. Now, from \cite[Formula (5.18)]{AmGrPa}
we know that
\begin{equation}
\left|D\left(\vert D_{x_{j}}u_{\varepsilon}\vert^{\frac{p-2}{2}}(D_{x_{j}}u_{\varepsilon})\left[\Phi(P)\right]^{\frac{p}{2}}\right)\right|^{2}\leq\,c_{6}(n,p)\,\left(\mathbf{A}_{1}+\mathbf{A}_{2}\right),\label{eq:A1+A2}
\end{equation}
where 
\[
\mathbf{A}_{1}:=\,\vert Du_{\varepsilon}\vert^{p-2}\,\vert D^{2}u_{\varepsilon}\vert^{2}\,[\Phi(P)]^{p}
\]
and 
\[
\mathbf{A}_{2}:=\,\vert Du_{\varepsilon}\vert^{p}\,\vert D^{2}u_{\varepsilon}\vert^{2}\,[\Phi(P)]^{p-2}\,[\Phi'(P)]^{2}.
\]
\foreignlanguage{english}{At this point, we estimate $\mathbf{A}_{1}$
and $\mathbf{A}_{2}$ separately in the set where $\vert Du_{\varepsilon}\vert>\lambda$,
since both $\mathbf{A}_{1}$ and $\mathbf{A}_{2}$ vanish in the set
$\{|Du_{\varepsilon}|\le\lambda\}$. Recalling the definitions of
$\Phi$ and $P$ in (\ref{eq:Phi}) and (\ref{eq:Pfunction}) respectively,
we can write $\mathbf{A}_{1}$ as the product of two terms: 
\[
\mathbf{A}_{1}=\,\frac{\vert Du_{\varepsilon}\vert^{p-1}\,(\vert Du_{\varepsilon}\vert-\lambda)_{+}^{(2\alpha-1)(p-1)}}{\left[\lambda^{2}+(\vert Du_{\varepsilon}\vert-\lambda)_{+}^{2}\right]^{\alpha(p-1)}}\,\cdot\,\frac{\vert D^{2}u_{\varepsilon}\vert^{2}\,(\vert Du_{\varepsilon}\vert-\lambda)_{+}^{p-1+2\alpha}}{\vert Du_{\varepsilon}\vert\left[\lambda^{2}+(\vert Du_{\varepsilon}\vert-\lambda)_{+}^{2}\right]^{\alpha}}\,.
\]
Then we have 
\[
\frac{\vert Du_{\varepsilon}\vert^{p-1}\,(\vert Du_{\varepsilon}\vert-\lambda)_{+}^{(2\alpha-1)(p-1)}}{\left[\lambda^{2}+(\vert Du_{\varepsilon}\vert-\lambda)_{+}^{2}\right]^{\alpha(p-1)}}\,\le\,2^{\alpha(p-1)}\,\frac{\vert Du_{\varepsilon}\vert^{2\alpha(p-1)}}{\vert Du_{\varepsilon}\vert^{2\alpha(p-1)}}\,=\,2^{\alpha(p-1)},
\]
where we have used that $2\alpha-1>0$, since by assumption $\alpha\ge\frac{p+1}{2(p-1)}>\frac{1}{2}$.
This implies that 
\begin{equation}
\mathbf{A}_{1}\,\le\,2^{\alpha(p-1)}\,\frac{\vert D^{2}u_{\varepsilon}\vert^{2}\,(\vert Du_{\varepsilon}\vert-\lambda)_{+}^{p-1+2\alpha}}{\vert Du_{\varepsilon}\vert\left[\lambda^{2}+(\vert Du_{\varepsilon}\vert-\lambda)_{+}^{2}\right]^{\alpha}}\,.\label{eq:A1}
\end{equation}
Now we deal with $\mathbf{A}_{2}$. First, we observe that 
\[
[\Phi(t)]^{p-2}\,[\Phi'(t)]^{2}=\frac{4\,\alpha^{2}\,\lambda^{4}\,t^{2(\alpha p-1)}}{(t^{2}+\lambda^{2})^{\alpha p+2}}\,.
\]
Therefore, $\mathbf{A}_{2}$ can be written as follows: 
\[
\mathbf{A}_{2}=\,4\,\alpha^{2}\,\lambda^{4}\,\frac{\vert Du_{\varepsilon}\vert^{p+1}\,(\vert Du_{\varepsilon}\vert-\lambda)_{+}^{2\alpha(p-1)-p-1}}{\left[\lambda^{2}+(\vert Du_{\varepsilon}\vert-\lambda)_{+}^{2}\right]^{\alpha(p-1)+2}}\,\cdot\,\frac{\vert D^{2}u_{\varepsilon}\vert^{2}\,(\vert Du_{\varepsilon}\vert-\lambda)_{+}^{p-1+2\alpha}}{\vert Du_{\varepsilon}\vert\left[\lambda^{2}+(\vert Du_{\varepsilon}\vert-\lambda)_{+}^{2}\right]^{\alpha}}\,.
\]
Again, the assumption $\alpha\ge\frac{p+1}{2(p-1)}$ implies $(\vert Du_{\varepsilon}\vert-\lambda)_{+}^{2\alpha(p-1)-p-1}\le\vert Du_{\varepsilon}\vert^{2\alpha(p-1)-p-1}$,
and so 
\[
\lambda^{4}\,\frac{\vert Du_{\varepsilon}\vert^{p+1}\,(\vert Du_{\varepsilon}\vert-\lambda)_{+}^{2\alpha(p-1)-p-1}}{\left[\lambda^{2}+(\vert Du_{\varepsilon}\vert-\lambda)_{+}^{2}\right]^{\alpha(p-1)+2}}\,\leq\,2^{\alpha(p-1)+2}\lambda^{4}\,\frac{\vert Du_{\varepsilon}\vert^{2\alpha(p-1)}}{\vert Du_{\varepsilon}\vert^{4}\,\vert Du_{\varepsilon}\vert^{2\alpha(p-1)}}\,<\,2^{\alpha(p-1)+2},
\]
where we have used that $\vert Du_{\varepsilon}\vert>\lambda$. Thus
we have 
\begin{equation}
\mathbf{A}_{2}\,\leq\,2^{\alpha(p-1)+4}\,\alpha^{2}\,\frac{\vert D^{2}u_{\varepsilon}\vert^{2}\,(\vert Du_{\varepsilon}\vert-\lambda)_{+}^{p-1+2\alpha}}{\vert Du_{\varepsilon}\vert\left[\lambda^{2}+(\vert Du_{\varepsilon}\vert-\lambda)_{+}^{2}\right]^{\alpha}}\,.\label{eq:A2}
\end{equation}
}Joining estimates \eqref{eq:new5}$-$(\ref{eq:A2}), for almost
every $t\in(t_{0}-\rho^{2},t_{1})$ we find\begin{align*}
[(D_{x_{j}}u_{\varepsilon}(\cdot,t))\,\eta^{2}\,\chi(t)\,\Phi(P(\cdot,t))]_{B_{p,\infty}^{2/p}(\mathbb{R}^{n})}^{p}\,&\leq\,c_{7}\int_{B_{\rho}}\frac{(\vert Du_{\varepsilon}(x,t)\vert-\lambda)_{+}^{p-1+2\alpha}\,\vert D^{2}u_{\varepsilon}(x,t)\vert^{2}}{\vert Du_{\varepsilon}(x,t)\vert\left[\lambda^{2}+(\vert Du_{\varepsilon}(x,t)\vert-\lambda)_{+}^{2}\right]^{\alpha}}\,\chi(t)\,\eta^{2}\,dx\\
&\,\,\,\,\,\,\,+\,\frac{c_{7}}{\rho^{2}}\,\chi(t)\int_{B_{\rho}}\vert Du_{\varepsilon}(x,t)\vert^{p}\,dx,
\end{align*} where $c_{7}\equiv c_{7}(n,p,\alpha)>0$. Combining the previous
inequality, (\ref{eq:extra}) and (\ref{eq:derf}) with \eqref{eq:pivot},
after some algebraic manipulation, from (\ref{eq:pivot00}) we get\begin{align*}
&\int_{Q^{t_{1}}}(D_{x_{j}}f_{\varepsilon})(D_{x_{j}}u_{\varepsilon})\,\eta^{2}(x)\,\chi(t)\,\Phi(P)\,dz\\
&\,\,\,\,\,\,\,\leq\,c_{7}\bigintsss_{t_{0}-\rho^{2}}^{t_{1}}\Vert f_{\varepsilon}(\cdot,t)\Vert_{B_{p',1}^{\frac{p-2}{p}}(B_{\rho})}\Biggl[\left(\int_{B_{\rho}}\frac{(\vert Du_{\varepsilon}(x,t)\vert-\lambda)_{+}^{p-1+2\alpha}\,\vert D^{2}u_{\varepsilon}(x,t)\vert^{2}}{\vert Du_{\varepsilon}(x,t)\vert\left[\lambda^{2}+(\vert Du_{\varepsilon}(x,t)\vert-\lambda)_{+}^{2}\right]^{\alpha}}\,\chi(t)\,\eta^{2}\,dx\right)^{\frac{1}{p}}\\
&\,\,\,\,\,\,\,\,\,\,\,\,\,\,\,\,\,\,\,\,\,\,\,\,\,\,\,\,\,\,\,\,\,\,\,\,\,\,\,\,\,\,\,\,\,\,\,\,\,\,\,\,\,\,\,\,\,\,\,\,\,\,\,\,\,\,\,\,\,\,\,\,\,\,\,\,\,\,\,\,\,\,\,\,\,+\left(1\,+\,\frac{1}{\rho^{2/p}}\right)\Vert Du_{\varepsilon}(\cdot,t)\Vert_{L^{p}(B_{\rho})}\Biggl]dt.
\end{align*}Now we go back to \eqref{eq:key1} and use the above estimate in combination
with Young's inequality. This yields\begin{align*}
&\int_{Q^{t_{1}}}\frac{(\vert Du_{\varepsilon}\vert-\lambda)_{+}^{p-1+2\alpha}}{\vert Du_{\varepsilon}\vert\,[\lambda^{2}+(\vert Du_{\varepsilon}\vert-\lambda)_{+}^{2}]^{\alpha}}\,\vert D^{2}u_{\varepsilon}\vert^{2}\,\chi(t)\,\eta^{2}(x)\,dz\\
&\,\,\,\,\,\,\,\leq\,c_{8}\left(1\,+\,\frac{1}{\rho^{2}}\right)\left[\int_{Q_{R}}(1+\vert Du_{\varepsilon}\vert^{2})^{\frac{p}{2}}\,dz\,+\,\Vert Du_{\varepsilon}\Vert_{L^{p}(Q_{R})}^{2}\right]+\,c_{8}\int_{t_{0}-\rho^{2}}^{t_{0}}\Vert f_{\varepsilon}(\cdot,t)\Vert_{B_{p',1}^{\frac{p-2}{p}}(B_{\rho})}^{p'}dt
\end{align*}for every $t_{1}\in(t_{0}-\rho^{2},t_{0})$ and a positive constant
$c_{8}$ depending only on $n$, $p$, $\alpha$ and $R$. Using the
properties of $\chi$ and $\eta$, from the above inequality we deduce\begin{align}\label{eq:new6}
&\int_{Q_{\rho/2}(z_{0})}\frac{(\vert Du_{\varepsilon}\vert-\lambda)_{+}^{p-1+2\alpha}}{\vert Du_{\varepsilon}\vert\,[\lambda^{2}+(\vert Du_{\varepsilon}\vert-\lambda)_{+}^{2}]^{\alpha}}\,\vert D^{2}u_{\varepsilon}\vert^{2}\,dz\nonumber\\
&\,\,\,\,\,\,\,\leq \left(c_{8}\,+\,\frac{c_{8}}{\rho^{2}}\right)\left[\int_{Q_{R}}(1+\vert Du_{\varepsilon}\vert^{2})^{\frac{p}{2}}\,dz\,+\,\Vert Du_{\varepsilon}\Vert_{L^{p}(Q_{R})}^{2}\right]+\,c_{8}\int_{t_{0}-\rho^{2}}^{t_{0}}\Vert f_{\varepsilon}(\cdot,t)\Vert_{B_{p',1}^{\frac{p-2}{p}}(B_{\rho})}^{p'}dt.
\end{align}At this point, recalling the definition of $\mathcal{V}_{\alpha,\lambda}$
in (\ref{eq:Vfun})$-$(\ref{eq:Gfun}), a straightforward computation
reveals that, for every $j\in\{1,\ldots,n\}$, we have\begin{align*}
D_{x_{j}}\mathcal{V}_{\alpha,\lambda}(Du_{\varepsilon})=&\,\,\,\frac{(\vert Du_{\varepsilon}\vert-\lambda)_{+}^{\frac{p-1+2\alpha}{2}}}{\vert Du_{\varepsilon}\vert^{2}\,[\lambda+(\vert Du_{\varepsilon}\vert-\lambda)_{+}]^{\frac{1+2\alpha}{2}}}\,\langle Du_{\varepsilon},D_{x_{j}}Du_{\varepsilon}\rangle\,Du_{\varepsilon}\\
&\,\,+\,\mathcal{G}_{\alpha,\lambda}((\vert Du_{\varepsilon}\vert-\lambda)_{+})\left[\frac{D_{x_{j}}Du_{\varepsilon}}{\vert Du_{\varepsilon}\vert}\,-\,\frac{\langle Du_{\varepsilon},D_{x_{j}}Du_{\varepsilon}\rangle}{\vert Du_{\varepsilon}\vert^{3}}\,Du_{\varepsilon}\right]
\end{align*}if $\vert Du_{\varepsilon}\vert>\lambda$, and $D_{x_{j}}\mathcal{V}_{\alpha,\lambda}(Du_{\varepsilon})=0$
otherwise. In the set $\{\vert Du_{\varepsilon}\vert>0\}$, this yields
\begin{equation}
\vert D_{x}\mathcal{V}_{\alpha,\lambda}(Du_{\varepsilon})\vert^{2}\leq\,\mathbf{B}_{1}+\mathbf{B}_{2}\,,\label{eq:B1+B2}
\end{equation}
where we define 
\[
\mathbf{B}_{1}:=\,2\,\frac{(\vert Du_{\varepsilon}\vert-\lambda)_{+}^{p-1+2\alpha}\,\vert D^{2}u_{\varepsilon}\vert^{2}}{[\lambda+(\vert Du_{\varepsilon}\vert-\lambda)_{+}]^{1+2\alpha}}
\]
and 
\[
\mathbf{B}_{2}:=\,8\,\frac{[\mathcal{G}_{\alpha,\lambda}((\vert Du_{\varepsilon}\vert-\lambda)_{+})]^{2}\,\vert D^{2}u_{\varepsilon}\vert^{2}}{\vert Du_{\varepsilon}\vert^{2}}\,.
\]
We now estimate $\mathbf{B}_{1}$ and $\mathbf{B}_{2}$ separately
in the set where $\vert Du_{\varepsilon}\vert>\lambda$, since both
$\mathbf{B}_{1}$ and $\mathbf{B}_{2}$ vanish for $0<\vert Du_{\varepsilon}\vert\leq\lambda$.
We immediately have
\begin{equation}
\mathbf{B}_{1}\leq\,2\,\frac{(\vert Du_{\varepsilon}\vert-\lambda)_{+}^{p-1+2\alpha}\,\vert D^{2}u_{\varepsilon}\vert^{2}}{\vert Du_{\varepsilon}\vert\,[\lambda^{2}+(\vert Du_{\varepsilon}\vert-\lambda)_{+}^{2}]^{\alpha}}\,.\label{eq:B1}
\end{equation}
As for $\mathbf{B}_{2}$, by Lemma \ref{lem:Glemma1} we obtain\begin{align}\label{eq:B2}
\mathbf{B}_{2}&\leq\,\frac{32}{p^{2}}\,\,\frac{(\vert Du_{\varepsilon}\vert-\lambda)_{+}^{p+1+2\alpha}\,\vert D^{2}u_{\varepsilon}\vert^{2}}{\vert Du_{\varepsilon}\vert^{2}\,[\lambda+(\vert Du_{\varepsilon}\vert-\lambda)_{+}]^{1+2\alpha}}\nonumber\\
&\le\,\frac{32}{p^{2}}\,\,\frac{(\vert Du_{\varepsilon}\vert-\lambda)_{+}^{p-1+2\alpha}\,(\vert Du_{\varepsilon}\vert-\lambda)_{+}^{2}\,\vert D^{2}u_{\varepsilon}\vert^{2}}{\vert Du_{\varepsilon}\vert^{3}\,[\lambda^{2}+(\vert Du_{\varepsilon}\vert-\lambda)_{+}^{2}]^{\alpha}}\nonumber\\
&\le\,\frac{32}{p^{2}}\,\,\frac{(\vert Du_{\varepsilon}\vert-\lambda)_{+}^{p-1+2\alpha}\,\vert D^{2}u_{\varepsilon}\vert^{2}}{\vert Du_{\varepsilon}\vert\,[\lambda^{2}+(\vert Du_{\varepsilon}\vert-\lambda)_{+}^{2}]^{\alpha}}\,.
\end{align}Joining estimates (\ref{eq:B1+B2})$-$\eqref{eq:B2}, we then find
\begin{equation}
\int_{Q_{\rho/2}(z_{0})}\vert D_{x}\mathcal{V}_{\alpha,\lambda}(Du_{\varepsilon})\vert^{2}\,dz\,\leq\,c(p)\int_{Q_{\rho/2}(z_{0})}\frac{(\vert Du_{\varepsilon}\vert-\lambda)_{+}^{p-1+2\alpha}\,\vert D^{2}u_{\varepsilon}\vert^{2}}{\vert Du_{\varepsilon}\vert\,[\lambda^{2}+(\vert Du_{\varepsilon}\vert-\lambda)_{+}^{2}]^{\alpha}}\,dz,\label{eq:combo}
\end{equation}
which combined with \eqref{eq:new6}, gives\begin{align*}
&\int_{Q_{\rho/2}(z_{0})}\vert D_{x}\mathcal{V}_{\alpha,\lambda}(Du_{\varepsilon})\vert^{2}\,dz\\
&\,\,\,\,\,\,\,\leq \left(c_{8}\,+\,\frac{c_{8}}{\rho^{2}}\right)\left[\int_{Q_{R}}(1+\vert Du_{\varepsilon}\vert^{2})^{\frac{p}{2}}\,dz\,+\,\Vert Du_{\varepsilon}\Vert_{L^{p}(Q_{R})}^{2}\right]+\,c_{8}\int_{t_{0}-\rho^{2}}^{t_{0}}\Vert f_{\varepsilon}(\cdot,t)\Vert_{B_{p',1}^{\frac{p-2}{p}}(B_{\rho})}^{p'}dt.
\end{align*}$\hspace*{1em}$Let us now consider the same $\varepsilon_{0}\in(0,1]$
as in Proposition \ref{prop:uniform1} and let $\varepsilon\in(0,\varepsilon_{0}]$.
Then, applying estimate (\ref{eq:uni1}), we obtain\begin{align*}
&\int_{Q_{\rho/2}(z_{0})}\vert D_{x}\mathcal{V}_{\alpha,\lambda}(Du_{\varepsilon})\vert^{2}\,dz\\
&\,\,\,\,\,\,\,\leq \left(C\,+\,\frac{C}{\rho^{2}}\right)\left[\Vert Du\Vert_{L^{p}(Q_{R})}^{p}+\Vert Du\Vert_{L^{p}(Q_{R})}^{2}+\lambda^{p}+\lambda^{2}+1\right]+\,C\,\Vert f_{\varepsilon}\Vert_{L^{p'}\left(t_{0}-\rho^{2},t_{0};B_{p',1}^{\frac{p-2}{p}}(B_{\rho})\right)}^{p'},
\end{align*} for a constant $C\equiv C(n,p,\alpha,R)>0$. Moreover, arguing as
in the proof of \cite[Proposition 5.1]{AmGrPa}, we find a positive
number $\varepsilon_{1}\leq\varepsilon_{0}$ such that\begin{align*}
\Vert f_{\varepsilon}\Vert_{L^{p'}\left(t_{0}-\rho^{2},t_{0};B_{p',1}^{\frac{p-2}{p}}(B_{\rho})\right)}^{p'}&=\int_{t_{0}-\rho^{2}}^{t_{0}}\Vert f_{\varepsilon}(\cdot,t)\Vert_{B_{p',1}^{\frac{p-2}{p}}(B_{\rho})}^{p'}dt\,\leq\,\int_{t_{0}-R^{2}}^{t_{0}}\Vert f(\cdot,t)\Vert_{B_{p',1}^{\frac{p-2}{p}}(B_{R})}^{p'}dt\\
&=\,\Vert f\Vert_{L^{p'}\left(t_{0}-R^{2},t_{0};B_{p',1}^{\frac{p-2}{p}}(B_{R})\right)}^{p'}<+\infty,\,\,\,\,\,\,\,\,\mathrm{for\,\,every\,\,}\varepsilon\in(0,\varepsilon_{1}].
\end{align*} Combining the last two estimates for $\varepsilon\in(0,\varepsilon_{1}]$,
we conclude the proof in the case $\lambda>0$.

\noindent $\hspace*{1em}$Finally, when $\lambda=0$ the above proof
can be greatly simplified, since the parameter $\alpha$ plays no
role and 
\[
\mathcal{V}_{0}(Du_{\varepsilon})=\,\frac{2}{p}\,H_{\frac{p}{2}}(Du_{\varepsilon})=\,\frac{2}{p}\,\vert Du_{\varepsilon}\vert^{\frac{p-2}{2}}Du_{\varepsilon}\,.
\]
In this regard, we leave the details to the reader.\end{proof}

\noindent $\hspace*{1em}$We conclude this section with the following
result, which will play a crucial role in the proof of Theorem \ref{thm:nuovo}.
\begin{prop}[\textbf{Uniform Sobolev estimate}]
\noindent  \label{prop:UniformSobolev-2} Under the assumptions of
Theorem \ref{thm:nuovo} and with the notation above, there exists
a positive number $\varepsilon_{0}\leq1$ such that, for every $\varepsilon\in(0,\varepsilon_{0}]$,
for every parabolic cylinder $Q_{\rho}(z_{0})\Subset Q_{R}(z_{0})$
and every $j\in\{1,\ldots,n\}$, we have
\begin{equation}
\int_{Q_{\rho/2}(z_{0})}\vert D_{x_{j}}\mathcal{V}_{\alpha,\lambda}(Du_{\varepsilon})\vert^{2}\,dz\,\leq\,\frac{C}{\rho^{2}}\left(\Vert Du\Vert_{L^{2}(Q_{R})}^{2}+\lambda^{2}+1\right)\,+\,C\left(\Vert f\Vert_{L^{2}(Q_{R})}^{2}+1\right)\label{eq:uniSob2}
\end{equation}
for a positive constant $C$ depending only on $n$ and $R$ in the
case $\lambda=0$, and additionally on $\alpha$ if $\lambda>0$. 
\end{prop}

\noindent \begin{proof}[\bfseries{Proof}] In what follows, we shall
keep the notation, definitions and choices used in the proof of Proposition
\ref{prop:UniformSobolev}.

\noindent $\hspace*{1em}$Let us first assume that $\lambda>0$. Letting
$\varrho\rightarrow0$ in \eqref{eq:new1} and arguing exactly as
in the preceding proof, we arrive at the following inequality
\begin{equation}
I_{3}\leq\,I_{4}+I_{5}+\tilde{I}_{6}+I_{7}+I_{8}\,,\label{eq:integrals-2}
\end{equation}
where\begin{align*}
&I_{3}:=\int_{Q^{t_{1}}}\sum_{j=1}^{n}\,\langle D^{2}A_{\varepsilon}(Du_{\varepsilon})\,D(D_{x_{j}}u_{\varepsilon}),D(D_{x_{j}}u_{\varepsilon})\rangle\,\chi(t)\,\eta^{2}(x)\,\Phi(P)\,dz,\\
&I_{4}:=\,4\int_{Q^{t_{1}}}\sum_{j=1}^{n}\,\langle D^{2}A_{\varepsilon}(Du_{\varepsilon})\,D\eta,D\eta\rangle\,\chi(t)\,\vert D_{x_{j}}u_{\varepsilon}\vert^{2}\,\Phi(P)\,dz,\\
&I_{5}:=\int_{Q^{t_{1}}}(\partial_{t}\chi)\,\eta^{2}(x)\left[\int_{0}^{\vert Du_{\varepsilon}(x,t)\vert^{2}}\Phi((\sqrt{w}-\lambda)_{+})\,dw\right]dz,\\
&\tilde{I}_{6}:=\,-\,4\int_{Q^{t_{1}}}\sum_{j=1}^{n}f_{\varepsilon}\,(D_{x_{j}}u_{\varepsilon})\,\eta(x)\,(D_{x_{j}}\eta)\,\chi(t)\,\Phi(P)\,dz,\\
&I_{7}:=\,-2\int_{Q^{t_{1}}}\sum_{j=1}^{n}f_{\varepsilon}\,(D_{x_{j}}^{2}u_{\varepsilon})\,\eta^{2}(x)\,\chi(t)\,\Phi(P)\,dz,\\
&I_{8}:=\,-2\int_{Q^{t_{1}}}\sum_{j=1}^{n}f_{\varepsilon}\,(D_{x_{j}}u_{\varepsilon})\,\eta^{2}(x)\,\chi(t)\,\Phi'(P)\,(D_{x_{j}}P)\,dz,
\end{align*}while $P$, $\eta$, $\chi$ and $\Phi$ are chosen as in the previous
proof. Note that we have already estimated $I_{3}$, $I_{4}$ and
$I_{5}$ in (\ref{eq:I3}), (\ref{eq:I4}) and (\ref{eq:I5}) respectively.
More precisely, for $p=2$ these estimates yield
\begin{equation}
I_{3}\,\geq\int_{Q^{t_{1}}}\frac{(\vert Du_{\varepsilon}\vert-\lambda)_{+}^{1+2\alpha}}{\vert Du_{\varepsilon}\vert\,[\lambda^{2}+(\vert Du_{\varepsilon}\vert-\lambda)_{+}^{2}]^{\alpha}}\,\vert D^{2}u_{\varepsilon}\vert^{2}\,\chi(t)\,\eta^{2}(x)\,dz,\label{eq:I3-2}
\end{equation}
\begin{equation}
I_{4}\,\leq\,\frac{8\,\tilde{c}^{2}}{\rho^{2}}\int_{Q_{R}}(1+\vert Du_{\varepsilon}\vert^{2})\,dz\label{eq:I4-2}
\end{equation}
and 
\begin{equation}
I_{5}\,\leq\,\frac{4}{3\rho^{2}}\,\Vert Du_{\varepsilon}\Vert_{L^{2}(Q_{R})}^{2}\,.\label{eq:I5-2}
\end{equation}
We now estimate $\tilde{I}_{6}$, $I_{7}$ and $I_{8}$ separately.
Let us first consider $\tilde{I}_{6}$. Using the fact that $0\leq\chi,\Phi\leq1$,
the properties (\ref{eq:proeta}) of $\eta$ and Young's inequality,
we obtain
\begin{equation}
\tilde{I}_{6}\,\leq\,4n\int_{Q_{R}}\vert f_{\varepsilon}\vert\,\vert Du_{\varepsilon}\vert\,\vert D\eta\vert\,dz\,\leq\,2n\,\Vert f_{\varepsilon}\Vert_{L^{2}(Q_{R})}^{2}\,+\,\frac{2n\,\tilde{c}^{2}}{\rho^{2}}\,\Vert Du_{\varepsilon}\Vert_{L^{2}(Q_{R})}^{2}\,.\label{eq:I6-2}
\end{equation}
As for the term $I_{7}$, recalling (\ref{eq:Pfunction}) and (\ref{eq:Phi}),
we have\begin{align*}
I_{7}\,&\leq\,2n\int_{Q^{t_{1}}\cap\{\vert Du_{\varepsilon}\vert\,>\,\lambda\}}\vert f_{\varepsilon}\vert\,\vert D^{2}u_{\varepsilon}\vert\,\eta^{2}(x)\,\chi(t)\,\frac{(\vert Du_{\varepsilon}\vert-\lambda)_{+}^{2\alpha}}{[\lambda^{2}+(\vert Du_{\varepsilon}\vert-\lambda)_{+}^{2}]^{\alpha}}\,dz\\
&\leq\,2n\int_{Q^{t_{1}}\cap\{\vert Du_{\varepsilon}\vert\,>\,\lambda\}}\vert f_{\varepsilon}\vert\,\vert D^{2}u_{\varepsilon}\vert\,\eta^{2}(x)\,\chi(t)\,\frac{(\vert Du_{\varepsilon}\vert-\lambda)_{+}^{2\alpha-1}\,\vert Du_{\varepsilon}\vert}{[\lambda^{2}+(\vert Du_{\varepsilon}\vert-\lambda)_{+}^{2}]^{\alpha}}\,dz.
\end{align*}Now we turn our attention to $I_{8}$. Recalling (\ref{eq:derPhi})
and applying Young's inequality, we get \begin{align*}
I_{8}\,&\leq\,4\alpha\int_{Q^{t_{1}}\cap\{\vert Du_{\varepsilon}\vert\,>\,\lambda\}}\sum_{j=1}^{n}\vert f_{\varepsilon}\vert\,\vert D_{x_{j}}u_{\varepsilon}\vert\,\eta^{2}(x)\,\chi(t)\,\frac{\lambda^{2}\,(\vert Du_{\varepsilon}\vert-\lambda)_{+}^{2\alpha-1}}{[\lambda^{2}+(\vert Du_{\varepsilon}\vert-\lambda)_{+}^{2}]^{\alpha+1}}\cdot\frac{\vert\langle Du_{\varepsilon},D_{x_{j}}Du_{\varepsilon}\rangle\vert}{\vert Du_{\varepsilon}\vert}\,dz\\
&\leq\,8\alpha n\int_{Q^{t_{1}}\cap\{\vert Du_{\varepsilon}\vert\,>\,\lambda\}}\vert f_{\varepsilon}\vert\,\vert D^{2}u_{\varepsilon}\vert\,\eta^{2}(x)\,\chi(t)\,\frac{\vert Du_{\varepsilon}\vert^{3}\,(\vert Du_{\varepsilon}\vert-\lambda)_{+}^{2\alpha-1}}{[\lambda^{2}+(\vert Du_{\varepsilon}\vert-\lambda)_{+}^{2}]^{\alpha}\,\vert Du_{\varepsilon}\vert^{2}}\,dz\\
&=\,8\alpha n\int_{Q^{t_{1}}\cap\{\vert Du_{\varepsilon}\vert\,>\,\lambda\}}\vert f_{\varepsilon}\vert\,\vert D^{2}u_{\varepsilon}\vert\,\eta^{2}(x)\,\chi(t)\,\frac{(\vert Du_{\varepsilon}\vert-\lambda)_{+}^{2\alpha-1}\,\vert Du_{\varepsilon}\vert}{[\lambda^{2}+(\vert Du_{\varepsilon}\vert-\lambda)_{+}^{2}]^{\alpha}}\,dz.
\end{align*}Joining the last two estimates and using Young's inequality again,
we find\begin{align}\label{eq:I7I8}
& I_{7}+I_{8}\nonumber\\
&\,\,\,\,\,\,\,\leq\,2n\,(4\alpha+1)\int_{Q^{t_{1}}\cap\{\vert Du_{\varepsilon}\vert\,>\,\lambda\}}\frac{\vert f_{\varepsilon}\vert\,\vert Du_{\varepsilon}\vert^{\frac{3}{2}}\,(\vert Du_{\varepsilon}\vert-\lambda)^{\alpha-\frac{3}{2}}}{[\lambda^{2}+(\vert Du_{\varepsilon}\vert-\lambda)^{2}]^{\frac{\alpha}{2}}}\cdot\frac{(\vert Du_{\varepsilon}\vert-\lambda)_{+}^{\alpha+\frac{1}{2}}\vert D^{2}u_{\varepsilon}\vert\,\eta^{2}(x)\,\chi(t)}{\sqrt{\vert Du_{\varepsilon}\vert}\,[\lambda^{2}+(\vert Du_{\varepsilon}\vert-\lambda)_{+}^{2}]^{\frac{\alpha}{2}}}\,dz\nonumber\\
&\,\,\,\,\,\,\,\leq\,2n^{2}(4\alpha+1)^{2}\int_{Q^{t_{1}}\cap\{\vert Du_{\varepsilon}\vert\,>\,\lambda\}}\frac{\vert f_{\varepsilon}\vert^{2}\,\vert Du_{\varepsilon}\vert^{3}\,(\vert Du_{\varepsilon}\vert-\lambda)^{2\alpha-3}\,\eta^{2}(x)\,\chi(t)}{[\lambda^{2}+(\vert Du_{\varepsilon}\vert-\lambda)^{2}]^{\alpha}}\,dz\nonumber\\
&\,\,\,\,\,\,\,\,\,\,\,\,\,\,+\,\frac{1}{2}\int_{Q^{t_{1}}}\frac{(\vert Du_{\varepsilon}\vert-\lambda)_{+}^{1+2\alpha}\,\vert D^{2}u_{\varepsilon}\vert^{2}}{\vert Du_{\varepsilon}\vert\,[\lambda^{2}+(\vert Du_{\varepsilon}\vert-\lambda)_{+}^{2}]^{\alpha}}\,\eta^{2}(x)\,\chi(t)\,dz\nonumber\\
&\,\,\,\,\,\,\,\leq\,4n^{2}(4\alpha+1)^{2}\int_{Q^{t_{1}}\cap\{\vert Du_{\varepsilon}\vert\,>\,\lambda\}}\frac{\vert f_{\varepsilon}\vert^{2}\,\vert Du_{\varepsilon}\vert^{3}\,(\vert Du_{\varepsilon}\vert-\lambda)_{+}^{2\alpha-3}}{\vert Du_{\varepsilon}\vert^{2\alpha}}\,dz\nonumber\\
&\,\,\,\,\,\,\,\,\,\,\,\,\,\,+\,\frac{1}{2}\int_{Q^{t_{1}}}\frac{(\vert Du_{\varepsilon}\vert-\lambda)_{+}^{1+2\alpha}\,\vert D^{2}u_{\varepsilon}\vert^{2}}{\vert Du_{\varepsilon}\vert\,[\lambda^{2}+(\vert Du_{\varepsilon}\vert-\lambda)_{+}^{2}]^{\alpha}}\,\eta^{2}(x)\,\chi(t)\,dz\nonumber\\
&\,\,\,\,\,\,\,\leq\,4n^{2}(4\alpha+1)^{2}\,\Vert f_{\varepsilon}\Vert_{L^{2}(Q_{R})}^{2}\,+\,\frac{1}{2}\int_{Q^{t_{1}}}\frac{(\vert Du_{\varepsilon}\vert-\lambda)_{+}^{1+2\alpha}\,\vert D^{2}u_{\varepsilon}\vert^{2}}{\vert Du_{\varepsilon}\vert\,[\lambda^{2}+(\vert Du_{\varepsilon}\vert-\lambda)_{+}^{2}]^{\alpha}}\,\eta^{2}(x)\,\chi(t)\,dz,
\end{align}where, in the last line, we have used the fact that $\vert Du_{\varepsilon}\vert^{3}\,(\vert Du_{\varepsilon}\vert-\lambda)_{+}^{2\alpha-3}\leq\vert Du_{\varepsilon}\vert^{2\alpha}$
in the set $\{\vert Du_{\varepsilon}\vert>\lambda\}$, since $\alpha\geq\frac{3}{2}$
by assumption. Now we combine estimates (\ref{eq:integrals-2})$-$\eqref{eq:I7I8},
thus obtaining
\[
\int_{Q^{t_{1}}}\frac{(\vert Du_{\varepsilon}\vert-\lambda)_{+}^{1+2\alpha}\,\vert D^{2}u_{\varepsilon}\vert^{2}}{\vert Du_{\varepsilon}\vert\,[\lambda^{2}+(\vert Du_{\varepsilon}\vert-\lambda)_{+}^{2}]^{\alpha}}\,\eta^{2}(x)\,\chi(t)\,dz\,\leq\,\frac{c_{1}}{\rho^{2}}\int_{Q_{R}}(1+\vert Du_{\varepsilon}\vert^{2})\,dz\,+\,c_{1}\,\Vert f_{\varepsilon}\Vert_{L^{2}(Q_{R})}^{2}\,,
\]
which holds for every $t_{1}\in(t_{0}-\rho^{2},t_{0})$ and a positive
constant $c_{1}$ depending only on $n$ and $\alpha$. Using the
properties of $\eta$ and $\chi$, from the above inequality we deduce
\begin{equation}
\int_{Q_{\rho/2}(z_{0})}\frac{(\vert Du_{\varepsilon}\vert-\lambda)_{+}^{1+2\alpha}\,\vert D^{2}u_{\varepsilon}\vert^{2}}{\vert Du_{\varepsilon}\vert\,[\lambda^{2}+(\vert Du_{\varepsilon}\vert-\lambda)_{+}^{2}]^{\alpha}}\,dz\,\leq\,\frac{c_{1}}{\rho^{2}}\int_{Q_{R}}(1+\vert Du_{\varepsilon}\vert^{2})\,dz\,+\,c_{1}\,\Vert f_{\varepsilon}\Vert_{L^{2}(Q_{R})}^{2}\,.\label{eq:combo-2}
\end{equation}
Taking $p=2$ into (\ref{eq:combo}) and combining the resulting estimate
with (\ref{eq:combo-2}), we then find
\[
\int_{Q_{\rho/2}(z_{0})}\vert D_{x}\mathcal{V}_{\alpha,\lambda}(Du_{\varepsilon})\vert^{2}\,dz\,\leq\,\frac{c_{1}}{\rho^{2}}\int_{Q_{R}}(1+\vert Du_{\varepsilon}\vert^{2})\,dz\,+\,c_{1}\,\Vert f_{\varepsilon}\Vert_{L^{2}(Q_{R})}^{2}\,.
\]
Let us now consider the same $\varepsilon_{0}\in(0,1]$ as in Proposition
\ref{prop:uniform1} and let $\varepsilon\in(0,\varepsilon_{0}]$.
Then, applying estimate (\ref{eq:uni1}) with $p=2$, we get
\[
\int_{Q_{\rho/2}(z_{0})}\vert D_{x}\mathcal{V}_{\alpha,\lambda}(Du_{\varepsilon})\vert^{2}\,dz\,\leq\,\frac{C}{\rho^{2}}\left(\Vert Du\Vert_{L^{2}(Q_{R})}^{2}+\lambda^{2}+1\right)\,+\,C\,\Vert f_{\varepsilon}\Vert_{L^{2}(Q_{R})}^{2}
\]
for a constant $C\equiv C(n,\alpha,R)>0$. Moreover, for $p=2$, from
(\ref{eq:fconv}) it follows that
\[
\Vert f_{\varepsilon}\Vert_{L^{2}(Q_{R})}\,\leq\,\Vert f\Vert_{L^{2}(Q_{R})}+1.
\]
Combining the last two inequalities, we conclude the proof in the
case $\lambda>0$.\\
$\hspace*{1em}$Finally, when $\lambda=0$ the above proof can be
greatly simplified, since the parameter $\alpha$ plays no role and
\[
\mathcal{V}_{0}(Du_{\varepsilon})\,=\,Du_{\varepsilon}\,.
\]
We leave the details to the reader.\end{proof}

\section{Proofs of Theorems \ref{thm:theo1} and \ref{thm:nuovo} \label{sec:proof1}}

\noindent $\hspace*{1em}$In this section, we prove Theorems \ref{thm:theo1}
and \ref{thm:nuovo} by combining a standard comparison argument (see
e.g. \cite{AmPa,Duzaar,GPdN}) with the estimates from Propositions
\ref{prop:uniform1}, \ref{prop:comparison}, \ref{prop:UniformSobolev}
and \ref{prop:UniformSobolev-2}.

\noindent \begin{proof}[\bfseries{Proof of Theorem~\ref{thm:theo1}}]
In what follows, we will refer to the notation introduced in Section
\ref{sec:a priori}. For any fixed $\varepsilon\in(0,1]$, we let
\[
u_{\varepsilon}\in C^{0}\left([t_{0}-R^{2},t_{0}];L^{2}(B_{R}(x_{0}))\right)\cap L^{p}\left(t_{0}-R^{2},t_{0};W^{1,p}(B_{R}(x_{0}))\right)
\]
be the unique energy solution of the Cauchy-Dirichlet problem (\ref{eq:CAUCHYDIR}).
Moreover, by a slight abuse of notation, for $w\in L_{loc}^{1}(Q_{R}(z_{0}),\mathbb{R}^{k})$,
$j\in\{1,\ldots,n\}$ and $h\neq0$, we set (when $x+he_{j}\in B_{R}(x_{0})$)
\[
\tau_{j,h}w(x,t):=\,w(x+he_{j},t)-w(x,t),
\]
\[
\Delta_{j,h}w(x,t):=\,\frac{w(x+he_{j},t)-w(x,t)}{h},
\]
where $e_{j}$ is the unit vector in the direction $x_{j}$. Now we
fix arbitrary radii
\[
0<r<\rho<R
\]
and use the finite difference operator $\tau_{j,h}$ defined above,
for increments $h\in\mathbb{R}\setminus\{0\}$ such that $\vert h\vert<\frac{\rho-r}{4}$.\\
$\hspace*{1em}$Let us first assume that $\lambda>0$. Furthermore,
consider the same $\varepsilon_{1}\in(0,1]$ as in Proposition \ref{prop:UniformSobolev}
and let $\varepsilon\in(0,\varepsilon_{1}]$. In the following, we
will denote by $C$ a positive constant which neither depends on $h$
nor on $\varepsilon$. In order to obtain an estimate for the finite
difference $\tau_{j,h}\mathcal{V}_{\alpha,\lambda}(Du)$, we use the
following comparison argument:\begin{align*}
&\int_{Q_{r/2}(z_{0})}\left|\tau_{j,h}\mathcal{V}_{\alpha,\lambda}(Du)\right|^{2}dz\\
&\,\,\,\,\,\,\,\leq\,4\int_{Q_{r/2}(z_{0})}\left|\tau_{j,h}\mathcal{V}_{\alpha,\lambda}(Du_{\varepsilon})\right|^{2}dx\,dt\,+\,4\int_{Q_{r/2}(z_{0})}\left|\mathcal{V}_{\alpha,\lambda}(Du_{\varepsilon})-\mathcal{V}_{\alpha,\lambda}(Du)\right|^{2}dx\,dt\\
&\,\,\,\,\,\,\,\,\,\,\,\,\,\,+\,4\int_{Q_{r/2}(z_{0})}\left|\mathcal{V}_{\alpha,\lambda}(Du_{\varepsilon}(x+he_{j},t))-\mathcal{V}_{\alpha,\lambda}(Du(x+he_{j},t))\right|^{2}dx\,dt\\
&\,\,\,\,\,\,\,\leq\,4\int_{Q_{r/2}(z_{0})}\left|\tau_{j,h}\mathcal{V}_{\alpha,\lambda}(Du_{\varepsilon})\right|^{2}dz\,+\,8\int_{Q_{R}(z_{0})}\left|\mathcal{V}_{\alpha,\lambda}(Du_{\varepsilon})-\mathcal{V}_{\alpha,\lambda}(Du)\right|^{2}dz.
\end{align*}Combining the previous estimate with \eqref{eq:comp} and \eqref{eq:uniSobolev},
for every $j\in\{1,\ldots,n\}$ we get \begin{align}\label{eq:mainest1}
&\int_{Q_{r/2}}\left|\tau_{j,h}\mathcal{V}_{\alpha,\lambda}(Du)\right|^{2}dz\nonumber\\
&\,\,\,\,\,\,\,\leq\,\frac{C(\rho^{2}+1)}{\rho^{2}}\,\vert h\vert^{2}\left(\Vert Du\Vert_{L^{p}(Q_{R})}^{p}+\Vert Du\Vert_{L^{p}(Q_{R})}^{2}+\lambda^{p}+\lambda^{2}+1\right)\nonumber\\
&\,\,\,\,\,\,\,\,\,\,\,\,\,\,+\,C\,\vert h\vert^{2}\,\Vert f\Vert_{L^{p'}\left(t_{0}-R^{2},t_{0};B_{p',1}^{\frac{p-2}{p}}(B_{R}(x_{0}))\right)}^{p'}\,+\,C\,\varepsilon\,\left(\Vert Du\Vert_{L^{p}(Q_{R})}^{p}+1\right)\nonumber\\
&\,\,\,\,\,\,\,\,\,\,\,\,\,\,+\,C\,\Vert f-f_{\varepsilon}\Vert_{L^{p'}(Q_{R})}\,\left(\Vert Du\Vert_{L^{p}(Q_{R})}+\lambda+1\right),
\end{align}which holds for every sufficiently small $h\in\mathbb{R}\setminus\{0\}$
and a constant $C\equiv C(n,p,\alpha,R)>0$. Therefore, recalling
(\ref{eq:strong}) and letting $\varepsilon\searrow0$ in \eqref{eq:mainest1},
we obtain\begin{align}\label{eq:mainest2}
\int_{Q_{r/2}}\left|\Delta_{j,h}\mathcal{V}_{\alpha,\lambda}(Du)\right|^{2}dz\,&\leq \left(C\,+\,\frac{C}{\rho^{2}}\right)\left[\Vert Du\Vert_{L^{p}(Q_{R})}^{p}+\Vert Du\Vert_{L^{p}(Q_{R})}^{2}+\lambda^{p}+\lambda^{2}+1\right]\nonumber\\
&\,\,\,\,\,\,\,+\,C\,\Vert f\Vert_{L^{p'}\left(t_{0}-R^{2},t_{0};B_{p',1}^{\frac{p-2}{p}}(B_{R}(x_{0}))\right)}^{p'}.
\end{align}Since the above estimate holds for every $j\in\{1,\ldots,n\}$ and
every suitably small $h\neq0$, by Lemma \ref{lem:RappIncre} we may
conclude that 
\[
\mathcal{V}_{\alpha,\lambda}(Du)\in L_{loc}^{2}\left(0,T;W_{loc}^{1,2}(\Omega,\mathbb{R}^{n})\right).
\]
Moreover, letting $h\rightarrow0$ in the previous inequality, we
also obtain estimate \eqref{eq:apriori1}.\\
\textcolor{red}{$\hspace*{1em}$}Let us now assume that $\lambda=0$.
In this case, we have already observed that the constant $C$ is independent
of $\alpha$ (see Proposition \ref{prop:UniformSobolev}) and 
\[
\mathcal{V}_{0}(Du)=\,\frac{2}{p}\,\vert Du\vert^{\frac{p-2}{2}}Du\,.
\]
Combining this fact with Lemma \ref{lem:Lind} and estimate \eqref{eq:mainest2},
for every sufficiently small $h\neq0$ we find 
\[
\int_{Q_{r/2}}\frac{\vert\tau_{j,h}Du(x,t)\vert^{p}}{\vert h\vert^{2}}\,dx\,dt\,\leq\,C_{1}(p)\int_{Q_{r/2}}\left|\Delta_{j,h}\mathcal{V}_{0}(Du)\right|^{2}dz\,\leq M
\]
for some finite positive constant $M$ depending on
\[
n,\,\,\,\,\,\,p,\,\,\,\,\,\,R,\,\,\,\,\,\,\rho,\,\,\,\,\,\,\Vert Du\Vert_{L^{p}(Q_{R})},\,\,\,\,\,\,\Vert f\Vert_{L^{p'}\left(t_{0}-R^{2},t_{0};B_{p',1}^{\frac{p-2}{p}}(B_{R})\right)},
\]
but not on $h$. Note that the above estimate also holds for every
$j\in\{1,\ldots,n\}$. Therefore, using Proposition \ref{prop:Nikolskii}
with the choices $G=Du$, $q=p$ and $\theta=\frac{2}{p}$, as well
as a standard covering argument, we infer that 
\[
Du\in L_{loc}^{p}\left(0,T;W_{loc}^{\sigma,p}(\Omega,\mathbb{R}^{n})\right)\,\,\,\,\,\,\,\,\,\mathrm{for}\,\,\mathrm{all}\,\,\sigma\in\left(0,\frac{2}{p}\right).
\]
This completes the proof.\end{proof}

\noindent $\hspace*{1em}$We now conclude this section by giving the 

\noindent \begin{proof}[\bfseries{Proof of Theorem~\ref{thm:nuovo}}]
Actually, it is sufficient to proceed as in the proof of Theorem \ref{thm:theo1},
but using estimate (\ref{eq:uniSob2}) in place of \eqref{eq:uniSobolev}.
In particular, since $\mathcal{V}_{0}(Du)=Du$, in the case $\lambda=0$
we retrieve 
\[
u\in L_{loc}^{2}(0,T;W_{loc}^{2,2}(\Omega))\,,
\]
which is a well-established result for the non-homogeneous heat equation
(\ref{eq:heat}) when $f\in L_{loc}^{2}(\Omega_{T}$).\end{proof}

\section{The time derivative: proof of Theorem \ref{thm:theo2} \label{sec:proof2}}

\noindent $\hspace*{1em}$This section is devoted to the study of
the existence and regularity of the time derivative of the weak solutions
to equation (\ref{eq:p-Poisson}), under the assumptions of Theorem
\ref{thm:theo2}. Indeed, we are now in a position to prove the aforementioned
theorem.

\noindent \begin{proof}[\bfseries{Proof of Theorem~\ref{thm:theo2}}]
We shall keep the notation and the parabolic cylinders used for the
proof of Theorem \ref{thm:theo1}. Let us notice that 
\[
\vert Du\vert^{p-2}Du=\mathcal{H}(\mathcal{V}_{0}(Du)),
\]
where $\mathcal{\mathcal{H}}:\mathbb{R}^{n}\rightarrow\mathbb{R}^{n}$
is the function defined by
\[
\mathcal{\mathcal{H}}(\xi):=\left(\frac{p}{2}\right)^{\frac{2}{p'}}\vert\xi\vert^{\frac{p-2}{p}}\xi\,,
\]
which is locally Lipschitz continuous for $p>2$. Thus, the vector
field $\vert Du\vert^{p-2}Du$ is weakly differentiable with respect
to the $x$-variable by virtue of the chain rule in Sobolev spaces.
From the definitions of $\mathcal{H}$ and $\mathcal{V}_{0}$, it
follows that
\begin{equation}
\left|D_{\xi}\,\mathcal{\mathcal{H}}(\mathcal{V}_{0}(Du))\right|\,\leq\,c_{1}\,\vert Du\vert^{\frac{p-2}{2}}\label{eq:estF}
\end{equation}
for some positive constant $c_{1}\equiv c(n,p)$. Now, applying the
chain rule, the Cauchy-Schwarz inequality and estimate (\ref{eq:estF}),
we obtain \begin{equation}\label{eq:new001}
\begin{split}
\left|D_{x}\,\mathcal{H}(\mathcal{V}_{0}(Du))\right|^{p'}&\leq\,c_{2}(n,p)\,\vert D_{\xi}\,\mathcal{H}(\mathcal{V}_{0}(Du))\vert^{p'}\,\vert D_{x}\mathcal{V}_{0}(Du)\vert^{p'}\\
&\leq\,c_{3}\,\vert Du\vert^{\frac{(p-2)p'}{2}}\,\vert D_{x}\mathcal{V}_{0}(Du)\vert^{p'},
\end{split}
\end{equation}where $c_{3}\equiv c_{3}(n,p)>0$. Using \eqref{eq:new001}, Hölder's
inequality with exponents $\left(\frac{2(p-1)}{p-2},\frac{2}{p'}\right)$
and estimate \eqref{eq:apriori1} with $\lambda=0$, we get\begin{align}\label{eq:new002}
\left(\int_{Q_{r/2}}\left|D_{x}\,\mathcal{H}(\mathcal{V}_{0}(Du))\right|^{p'}dz\right)^{\frac{1}{p'}}& \leq\,c_{4}(n,p)\,\Vert Du\Vert_{L^{p}(Q_{r/2})}^{\frac{p-2}{2}}\,\left(\int_{Q_{r/2}}\vert D_{x}\mathcal{V}_{0}(Du)\vert^{2}\,dz\right)^{\frac{1}{2}}\nonumber\\
&\leq \left(C\,+\,\frac{C}{\rho}\right)\left[\Vert Du\Vert_{L^{p}(Q_{R})}^{p-1}+\Vert Du\Vert_{L^{p}(Q_{R})}^{\frac{p}{2}}+\Vert Du\Vert_{L^{p}(Q_{R})}^{\frac{p-2}{2}}\right]\nonumber\\
&\,\,\,\,\,\,\,+\,C\,\Vert Du\Vert_{L^{p}(Q_{R})}^{\frac{p-2}{2}}\,\Vert f\Vert_{L^{p'}\left(t_{0}-R^{2},t_{0};B_{p',1}^{\frac{p-2}{p}}(B_{R})\right)}^{\frac{p'}{2}}
\end{align}for a constant $C\equiv C(n,p,R)>0$. Note that the right-hand side
of \eqref{eq:new002} is finite, and this implies that 
\begin{equation}
\vert Du\vert^{p-2}Du\in L_{loc}^{p'}(0,T;W_{loc}^{1,p'}(\Omega,\mathbb{R}^{n})).\label{eq:higherdiff}
\end{equation}
Therefore, going back to the weak formulation (\ref{eq:weaksol}),
we can now perform a partial integration in the second term on the
left-hand side with respect to the spatial variables. We thus obtain
\[
\int_{Q_{r/2}(z_{0})}u\cdot\partial_{t}\varphi\,\,dz\,=\,-\int_{Q_{r/2}(z_{0})}\left(\sum_{j=1}^{n}D_{x_{j}}[\vert Du\vert^{p-2}D_{x_{j}}u]+f\right)\cdot\varphi\,\,dz
\]
for any $\varphi\in C_{0}^{\infty}(Q_{r/2}(z_{0}))$, and the desired
conclusion immediately follows from (\ref{eq:higherdiff}), since
\[
f\in L_{loc}^{p'}\left(0,T;B_{p',1,loc}^{\frac{p-2}{p}}(\Omega)\right).
\]
$\hspace*{1em}$Furthermore, we can now observe that
\[
\partial_{t}u\,=\,\sum_{j=1}^{n}D_{x_{j}}[\vert Du\vert^{p-2}D_{x_{j}}u]+f\,\,\,\,\,\,\,\,\,\mathrm{in}\,\,\,\,Q_{r/2}(z_{0}).
\]
From this we can infer
\[
\left(\int_{Q_{r/2}(z_{0})}\left|\partial_{t}u\right|^{p'}dz\right)^{\frac{1}{p'}}\leq\,\,n\,\Vert D_{x}\,\mathcal{H}(\mathcal{V}_{0}(Du))\Vert_{L^{p'}(Q_{r/2})}+\,\Vert f\Vert_{L^{p'}(Q_{r/2})}\,,
\]
and combining this inequality with \eqref{eq:new002}, we finally
obtain estimate \eqref{eq:stimaapriori2}. \end{proof}\medskip{}

\noindent \textbf{Acknowledgments. }The author has been partially
supported by the INdAM$-$GNAMPA 2024 Project “Fenomeno di Lavrentiev,
Bounded Slope Condition e regolarità per minimi di funzionali integrali
con crescite non standard e lagrangiane non uniformemente convesse”
(CUP: E53C23001670001).\bigskip{}

\noindent \textbf{Data availability.} There are no data associated
with this manuscript.\\
 \\
\textbf{Declarations.} The author states that there is no conflict
of interest.

\lyxaddress{\noindent \textbf{$\quad$}}

\begin{thebibliography}{10}
\bibitem{Adams}\foreignlanguage{english}{ R.A. Adams, \textit{Sobolev
Spaces}, Pure Appl. Math., vol. 65, Academic Press, New York, 1975.}

\bibitem{Amb1} P. Ambrosio, \textit{Besov regularity for a class
of singular or degenerate elliptic equations}, J. Math. Anal. Appl.
\textbf{505} (2022), no. 2, Paper No. 125636.

\bibitem{Amb2} P. Ambrosio, \textit{Fractional Sobolev regularity
for solutions to a strongly degenerate parabolic equation}, Forum
Math. \textbf{35} (6) (2023), 1485-1497.

\bibitem{AmBa} P. Ambrosio, F. Bäuerlein, \textit{Gradient bounds
for strongly singular or degenerate parabolic systems}, J. Differ.
Equ. \textbf{401} (2024), 492-549.

\bibitem{AmCuDe} P. Ambrosio, S. Cuomo, M. De Rosa, \textit{A physics-informed
deep learning approach for solving strongly degenerate parabolic problems},
Engineering with Computers, (2024). DOI: \url{https://doi.org/10.1007/s00366-024-01961-9}.

\bibitem{AmGrPa} P. Ambrosio, A.G. Grimaldi, A. Passarelli di Napoli,
\textit{Sharp second-order regularity for widely degenerate elliptic
equations}, preprint (2024), version 3. Available at \url{https://arxiv.org/abs/2401.13116v3}.

\bibitem{AmPa} P. Ambrosio, A. Passarelli di Napoli, \textit{Regularity
results for a class of widely degenerate parabolic equations}, Adv.
Calc. Var. \textbf{17} (3) (2024), 805-829.

\bibitem{BDGP0} \foreignlanguage{english}{V. Bögelein, F. Duzaar,
R. Giova, A. Passarelli di Napoli, \textit{Gradient regularity for
a class of widely degenerate parabolic systems}, SIAM J. Math. Anal.,
1-62 (2024).}

\bibitem{BDGP} \foreignlanguage{english}{V. Bögelein, F. Duzaar,
R. Giova, A. Passarelli di Napoli, \textit{Higher regularity in congested
traffic dynamics}, Math. Ann. \textbf{385}, 1-56 (2023).}

\bibitem{BrLIP} L. Brasco, \textit{Global $L^{\infty}$ gradient
estimates for solutions to a certain degenerate elliptic equation},
Nonlinear Anal., Theory Methods Appl. \textbf{74} (2) (2011) 516-531.

\bibitem{Br} L. Brasco, G. Carlier, F. Santambrogio, \textit{Congested
traffic dynamics, weak flows and very degenerate elliptic equations},
{[}corrected version of mr2584740{]}, J. Math. Pures Appl. (9) \textbf{93}
(2010), no. 6, 652-671.

\bibitem{BraSan} \foreignlanguage{english}{L. Brasco, F. Santambrogio,
\textit{A sharp estimate à la Calderón-Zygmund for the $p$-Laplacian},
Commun. Contemp. Math. \textbf{20} (2018), no. 3, Article ID 1750030.}

\bibitem{DuGaMi} F. Duzaar, A. Gastel, G. Mingione, \textit{Elliptic
systems, singular sets and Dini continuity}, Comm. Part. Diff. Equ.,\textbf{
29} (2004), 1215-1240.

\bibitem{Duzaar} F. Duzaar, G. Mingione, K. Steffen, \textit{Parabolic
systems with polynomial growth and regularity}, Mem. Amer. Math. Soc.
214 (2011).

\bibitem{EspLeoMin} L. Esposito, F. Leonetti, G. Mingione, \textit{Higher
integrability for minimizers of integral functionals with $(p,q)$
growth}, J. Differ. Equ. \textbf{157 }(1999), 414-438.

\bibitem{GPdN} A. Gentile, A. Passarelli di Napoli, \textit{Higher
regularity for weak solutions to degenerate parabolic problems}, Calc.
Var. \textbf{62}, 225 (2023). 

\bibitem{GiPaSc} F. Giannetti, A. Passarelli di Napoli, C. Scheven,
\textit{Higher differentiability of solutions of parabolic systems
with discontinuous coefficients}, J. Lond. Math. Soc. (2) \textbf{94}
(2016), no. 1, 1-20.

\bibitem{Giu} E. Giusti, \textit{Direct Methods in the Calculus of
Variations}, World Scientific Publishing Co., 2003.

\bibitem{Lind1} P. Lindqvist, \textit{On the time derivative in a
quasilinear equation}, Skr. K. Nor. Vidensk. Selsk. (2008), no. 2,
1-7.

\bibitem{Lind2} P. Lindqvist, \textit{On the time derivative in an
obstacle problem}, Rev. Mat. Iberoam. \textbf{28} (2012), no. 2, 577-590. 

\bibitem{Lind3} P. Lindqvist, \textit{The time derivative in a singular
parabolic equation}, Differential Integral Equations \textbf{30} (2017),
no. 9-10, 795-808.

\bibitem{Lind} P. Lindqvist, \textit{Notes on the Stationary $p$-Laplace
Equation}, Springer Briefs Math., Springer, Cham 2019.

\bibitem{Lions} J.-L. Lions, \textit{Quelques méthodes de résolution
des problèmes aux limites non linéaires}, Gauthier-Villars, Paris,
1969.

\bibitem{Ming1} G. Mingione, \textit{Gradient potential estimates},
J. European Math. Soc. \textbf{13} (2011), 459-486. 

\bibitem{Ming2} G. Mingione, \textit{The Calderón-Zygmund theory
for elliptic problems with measure data}, Ann. Scu. Norm. Sup. Pisa
Cl. Sci. (V) \textbf{6} (2007), 195-261.

\bibitem{Ming3} G. Mingione, \textit{The singular set of solutions
to non-differentiable elliptic systems}, Arch. Rat. Mech. Anal. \textbf{166
}(2003), 287-301.

\bibitem{Sche} C. Scheven, \textit{Regularity for subquadratic parabolic
systems: Higher integrability and dimension estimates}, Proc. Roy.
Soc. Edinburgh Sect. A \textbf{140} (2010), no. 6, 1269-1308.

\bibitem{Tri0}\foreignlanguage{english}{ H. Triebel, \textit{Spaces
of Besov-Hardy-Sobolev type}, Teubner-Texte Math. 15, Teubner, Leipzig,
1978.}

\bibitem{Triebel}\foreignlanguage{english}{ H. Triebel, \textit{Theory
of Function Spaces}, Monogr. Math. 78, Birkhäuser, Basel, 1983.}

\bibitem{TriebelIV}\foreignlanguage{english}{ H. Triebel, \textit{Theory
of Function Spaces IV}, Monogr. Math. 107, Birkhäuser, Basel, 2020.}
\end{thebibliography}
\end{document}